%% file: MR-Lasso.tex
\documentclass[12pt]{article} 

\usepackage[dvips]{graphics}
\DeclareGraphicsExtensions{.eps.gz,.eps,.epsi.gz,.epsi,.ps,.ps.gz}
\usepackage{epsfig}
\usepackage{color}
\graphicspath{{fig/}} \textwidth=6.5in \textheight=9in

\usepackage{natbib}
\usepackage[colorlinks,allcolors=blue]{hyperref}
\usepackage{latexsym}
\usepackage{graphicx}
\usepackage{amsmath}
\usepackage{amsfonts}
\usepackage{amssymb}
\usepackage{mathrsfs}
\usepackage{mathtools}
\usepackage[title,titletoc,toc]{appendix}
\usepackage{algorithm}
\usepackage{algorithmic}
\usepackage{hyperref}

\newtheorem{theorem}{Theorem}
\newtheorem{lemma}{Lemma}
\newtheorem{proposition}{Proposition}
\newtheorem{corollary}{Corollary}
\newtheorem{remark}{Remark}

\usepackage[dvipsnames]{xcolor}

\input{04-pre-am-short.tex}

\def\fbar{{\overline f}}

\def\RE{\hbox{\rm RE}}

\def\bffbar{{\overline \bff}}

\def\BR{{\rm\tiny BR}}

\def\REbar{{\overline \RE}}

\begin{document}

\title{Adaptive Estimation\\ In High-Dimensional Additive Models\\
With Multi-Resolution Group Lasso}
  
\author{Yisha Yao and 
 Cun-Hui Zhang\thanks{Partially supported by NSF Grants DMS-1721495, IIS-1741390 and CCF-1934924} 
   \ \\
Rutgers University\\ \ } 

\date{}

\maketitle

\begin{abstract}
In additive models with many nonparametric components, 
a number of regularized estimators have been proposed and proven to attain various  
error bounds under different combinations of sparsity and fixed smoothness conditions.  
Some of these error bounds match minimax rates in the corresponding settings. 
Some of the rate minimax methods are non-convex and 
computationally costly. 
From these perspectives, the existing solutions to the high-dimensional 
additive nonparametric regression problem are fragmented. 
In this paper, we propose a multi-resolution group Lasso (MR-GL) method in a unified approach 
to simultaneously achieve or improve existing error bounds and provide new ones 
without the knowledge of the level of sparsity or the degree of smoothness of the unknown functions. 
Such adaptive convergence rates are established when a prediction factor can be treated as a constant. 
Furthermore, we prove that the prediction factor, which can be bounded in terms of a restricted eigenvalue or a compatibility coefficient, can be indeed treated as a constant 
for random designs under a nearly optimal sample size condition. 
\end{abstract}


\section{Introduction}

Additive model (AM) dates back to the 1980's 
when the curse of dimensionality was a major concern in nonparametric regression, 
as the multivariate nonparametric model generally fails when there are not enough observations to fit a moderately smooth but otherwise unrestricted function of multiple inputs. 
The AM \citep{friedman1981projection,stone1985additive,hastie1986generalized} can be written as 
\bel{LM-0}
y_i = \sum_{j=1}^p f_j(x_{i,j}) + \veps_i,\quad  i=1,\ldots,n, 
\eel
where $\veps_i$ is the noise, $x_{i,j}$ is the $j$-th design variable with the $i$-th sample point,  
and each $f_j(\cdot)$ is an unknown univariate function satisfying a certain smoothness condition. 
Hereafter, we will call each $f_j$ a component function and their sum the regression function. 

For fixed $p$, \cite{stone1985additive} established a squared error rate $n^{-2\alpha/(2\alpha+1)}$ for the spline estimates of the regression function in (\ref{LM-0}), where $\alpha$ is the assumed degree of smoothness of the underlying functional class.
During the same period, numerous advances emerged in the study of specific 
additive non- and semi-parametric models and their extensions  
\citep{engle1986semiparametric, wahba1986partial, 
hastie1987generalized}. 
Still, model (\ref{LM-0}) remains among the most widely used and intensively studied AM owing to its balance of flexibility, interpretability, and tractability. 
When accounting for interaction terms, (\ref{LM-0}) 
leads to the smoothing spline analysis of variance (SS-ANOVA) model \citep{wahba1995smoothing}. 
When certain smoothness conditions are imposed on each component function, 
(\ref{LM-0}) is approximately the sum of several univariate smooth splines \citep{wahba1990spline, wood2004stable}. 
Along with the theoretical explorations 
of model (\ref{LM-0}), several algorithms have been proposed to solve the nonparametric AM problem 
including the alternating conditional expectations 
\citep{breiman1985estimating} and backfitting \citep{buja1989linear}. 
These algorithms construct a general iterative framework 
where in each iteration the component functions are estimated by local smoothing \citep{hastie1986generalized}, maximum likelihood \citep{hastie1987generalized} and other methods. 
When the component functions are restricted to smoothing splines, one may use generalized cross validation 
and generalized maximum likelihood 
\citep{craven1978smoothing, wahba1985comparison, gu1990adaptive, gu1991minimizing}. 
Iteratively re-weighted least squares 
\citep{green1984iteratively} and restricted maximum likelihood 
\citep{wood2011fast} are available when the AM is interpreted 
as a generalized linear model. 

Upon the emergence of the SS-ANOVA, 
researchers started to concern about the problem of component selection, removing 
the component functions that explain a little fraction of the model. The early 
approaches include the angle cosine between the residue vector and the components \citep{gu1992diagnostics}, Bayesian prior on the coefficients \citep{yau2003bayesian} or the component functions \citep{wood2002model}, sparse kernels \citep{gunn2002structural} and likelihood basis pursuit \citep{zhang2004variable}. 
Despite the reported efficiency in performance, 
these approaches are either hard to implement or short of theoretical guarantees. 
Component selection is also the main theme in the study of high-dimensional sparse additive models where  the sample size $n$ is much smaller than the number of additive components $p$ but much larger than the number of non-zero components. To seek a parsimonious model with better interpretation, 
component functions with little explanatory contribution should be eliminated. 

High-dimensional AM has been very active topic. Several different approaches have been introduced in various contexts, providing distinct interpretations. \cite{lin2006component} proposed  the COmponent Selection and Smoothing Operator (COSSO) in the context of reproducing kernel Hilbert space (RKHS) which essentially minimizes a penalized empirical risk function with the penalty term being proportional to the sum of the RKHS norms. The COSSO has the interpretations as a smoothing spline regularized by nonnegative garrote \citep{breiman1995better} as well as a generalization of the Lasso, and can be computed by iterations between fitting the smoothing spline and solving the nonnegative garrote. 
Both \cite{yuan2007nonnegative} and \cite{ravikumar2009sparse} proposed to apply a generalization of the 
nonnegative garrote to the sparse AM with penalization schemes analogous to that in \cite{lin2006component}, but different algorithms were implemented. 
In an selection consistency analysis, \cite{ravikumar2009sparse} interpreted their scheme as a functional version of the group Lasso after transforming each component function in basis expansion.  
\cite{huang2010variable} applied group Lasso to basis expansions of the components $f_j$ and provided sufficient conditions for selection consistency. 

Among more recent studies of the sparse AM, 
this paper is most closely related to the following papers. 
To abbreviate the discussion, assume here that the AM \eqref{LM-0} holds for $f_j(\cdot)=f^*_j(\cdot)$ 
and iid $\veps_i\sim N(0,\sigma^2)$ and that $f^*_j(\cdot)$ are uniformly $\alpha$-smooth in the sense 
that the $L_2$ norms of the $\alpha$-th derivative of $f^*_j(\cdot)$ (or the RKHS norms) are uniformly bounded, 
along with additional side conditions imposed in the respective following referenced papers. 
Let $f^* = f^*(x_1,\ldots,x_p)=\sum_{j=1}^p f^*_j(x_j)$ be the true regression function 
and $s_0=\#\{j: f^*_j(\cdot)\neq 0\}$. 
\cite{meier2009high} proved that their penalized least squares estimator (LSE) of $f^*$ achieves the squared error rate $s_0\{(\log p)/n\}^{2\alpha/(2\alpha+1)}$ under a compatibility condition, and for certain random designs the compatibility condition holds 
when $s_0\{(\log p)/n\}^{2\alpha/(2\alpha+1)-1/2}$ is sufficiently small. 
\cite{koltchinskii2010sparsity} established the squared error rate 
\bel{KY}
s_0\big\{n^{-2\alpha/(2\alpha+1)} + (\log p)/n\big\}
\eel
for their penalized LSE of $f^*$ 
under an $\ell_\infty$ constraint on $f^*$ and its estimator. 
As $n^{-2\alpha/(2\alpha+1)}$ is the minimax rate for the estimation of a single $\alpha$-smooth 
$f_j$ and $s_0(\log p)/n$ is the minimax rate when $f_j$ are all linear, 
\eqref{KY} is the minimax rate in the sparse AM without $\ell_\infty$ constraint. 
\cite{raskutti2012minimax} proved that the $\ell_\infty$ constraint in \cite{koltchinskii2010sparsity} 
can be removed when the covariates are independent uniform variables on $[0,1]$, 
and pointed out that when $\|f^*\|_\infty$ is 
bounded the squared error rate 
\bel{RWY}
(\log s_0)^{1/2}(s_0^{1/\alpha}/n)^{1/2} + (s_0/n)\log(p/s_0)
\eel 
is also achievable so that 
(\ref{KY}) is not necessarily sharp. 
\cite{suzuki2012fast} extended and improved the results of \cite{koltchinskii2010sparsity} 
under $\ell_1$ and $\ell_2$ conditions on the RKHS norms of $f^*_j(\cdot)$ and an $\ell_\infty$ 
condition on the noise $\veps_i$ in \eqref{LM-0}. 
\cite{yuan2015minimax} considered the LSE under a bounded $\ell_q$ constraint 
on the RKHS norms of $f^*_j(\cdot)$, $q\in (0,1)$, and proved squared error rate 
\bel{YZ}
n^{-2\alpha/(2\alpha+1)} + \big\{(\log p)/n\big\}^{1-q/2},
\eel
which is rate minimax for the $\ell_q$ class. 
\cite{tan2019doubly} considered regularized LSE with a sum of empirical and functional penalties 
and proved that the estimator achieves the minimax squared error rate 
\bel{TZ}
   n^{-(2-q)/(2+(1-q)/\alpha)} + \big\{(\log p)/n\big\}^{1-q/2}
\eel
under a bounded $\ell_1$ constraint on the RKHS norms of $f^*_j(\cdot)$ 
and a bounded $\ell_q$ constraint on the $L_2$ norms of $f^*_j(\cdot)$. 
In this literature, the penalty levels are chosen according to the function class under consideration, 
and the theoretical results all require additional side conditions. 

While the above results are very interesting in and of themselves individually and complementary to each other as a whole, the state of art solution to the high-dimensional additive nonparametric regression problem is still fragmented. 
A rate minimax estimator under a certain complexity condition on the true $f^*$ 
would lose to another estimator under a different 
complexity condition. 
Moreover, these existing methods all require \textit{a priori} information on the degree of smoothness of the underlying component functions $f^*_j(\cdot)$. 
Yet in practice, we often do not possess such information. 

In this paper, we propose a new method, Multi-Resolution Group Lasso (MR-GL), as an adaptive solution to the high-dimensional additive nonparametric regression problem. The MR-GL estimator will be shown to simultaneously achieve or improve existing error bounds, and it requires no knowledge of either the level of sparsity or the degree of smoothness. Thus, the method is adaptive and expected to be superior over existing approaches when the sparsity and the degree of smoothness of the true regression function are unknown. Furthermore, the MR-GL is easy to compute as it can be directly 
implemented as a group Lasso algorithm. 

In nonparametric regression, i.e. $p=1$ in \eqref{LM-0}, adaptive estimation dates back to 
\cite{efroimovich1984learning, efroimovich1986self} who proposed to apply shrinkage estimators 
to coefficient blocks in an orthonormal basis expansion of the unknown 
and proved the adaptive optimality of the method 
in the sense of simultaneous approximate minimaxity 
over a broad spectrum of Sobolev classes of regression functions. 
Such adaptive optimality can be also achieved with different estimators in more general classes, 
for example, block James-Stein estimators in Sobolev classes 
\citep{brown1997superefficiency,cai1999adaptive}, 
threshold estimators for near or rate minimaxity in more general Besov classes 
\citep{donoho1994ideal,donoho1998minimax,hall1999minimax}, and block general empirical Bayes 
estimators for exact minimaxity in Besov classes \citep{zhang2005general}. 
In Fourier basis such blocks represent disjoint bands of frequency, whereas
in wavelet bases such blocks represent different scales. 
Thus, the above adaptive estimation methods can be all viewed as 
multi-resolution \citep{mallat1989multiresolution,meyer1992wavelets} solutions to 
the nonparametric regression problem.

The thrust of this paper is that adaptive rate optimal estimation can be achieved in high-dimensional AM 
by directly applying the group Lasso in a carefully architected multi-resolution analysis, i.e.\,MR-GL. 
We write each nonparametric component $f_j(\cdot)$ as an infinite sum 
$f_j(x) = \sum_{k=k_*}^\infty f_{j,k}(x)$ with 
$f_{j,k}(x)= \sum_{\ell=1}^{2^{(k-1)\vee k_*}} u_{j,k,\ell}(x)\beta_{j,k,\ell}$ 
where $u_{j,k,\ell}(\cdot)$ are certain basis functions for $f_j(\cdot)$. 
After the basis expansion, we fit the response vector $\by = (y_1,\ldots,y_n)^\top$ 
by the group Lasso with $\bff_{j,k}=(f_{j,k}(x_{1,j}),\ldots,f_{j,k}(x_{n,j}))^\top$ 
as the group effect representing $f_j(\cdot)$ at the $k$-th resolution level (or in the $k$-th block of the basis), 
$k=k_*,\ldots,k^*$ and $j=1,\ldots,p$. We choose $k_*$ and $k^*$ to cover a wide 
spectrum of resolution levels so that adaptation is achieved for nonparametric components 
with heterogeneous smoothness indices $\alpha_j >1/2$. 

The group Lasso is not new in the analysis of the AM as we have briefly reviewed. 
As the group penalty can be viewed as a sum of inner-product norms of groups of coefficients in a linear system, most such penalized MLE approaches to the sparse AM use a group Lasso penalty on $f_j$ or a sum of two such penalties. 
However, as shown in the nonparametric regression literature and the sparse AM literature, 
regularizing one or two inner product norms of $f_j$ amounts to fixing the bandwidth for its estimation 
and thus are non-adaptive to its level of smoothness. 
Our analysis reveals that by applying group regularization to $f_{j,k}$ in a suitable multi-resolution 
structure, the optimality of the group Lasso translates to adaptive optimality in the estimation of the 
regression function in high-dimensional sparse AM.   


%

\medskip
{\bf Notation:}  
We write the data as $(\bX,\by) = (\bx_1,\ldots,\bx_p,\by)$ with covariate vectors $\bx_j=(x_{1,j},\ldots,x_{n,j})^\top$ as columns of $\bX$ and response vector $\by=(y_1,\ldots,y_n)^\top$. 
Moreover, we use bold face $\bg_j$ to denote realizations of function $g(x)$ in the sample vector $\bx_j$  
so that $\bg_j = g(\bx_j) = (g(x_{1,j}),\ldots, g(x_{n,j}))^\top$. Similarly for any function $g: \R^p\to\R$, 
$\bg = g(\bX)$ denotes the realization of $g$ in the sample as the vector with elements 
$g(x_{i,1},\ldots,x_{i,p})$. 
Thus, the vector version of the AM \eqref{LM-0} is $\bff = \sum_{j=1}^p \bff_j$ and its multi-resolution approximation can be written as $\sum_{j=1}^p\sum_{k=k_*}^{k^*}\bff_{j,k}$.
For vectors $\bu$ and $\bv\in \R^n$, $\langle \bu, \bv\rangle_n = \bu^\top\bv/n$ and 
$\|\bu\|_{2,n} = \|\bu\|_2/\sqrt{n}$ denote the length normalized $\ell_2$ inner product and norm. 
For random design $\bX$ and functions $g(\cdot)$ and $h(\cdot)$ from $\R^p$ to $\R$, 
$\langle g, h\rangle_{L_2} = \E\big[\langle \bg, \bh\rangle_n\big]$ and 
$\|g\|_{L_2} = \langle g, g\rangle_{L_2}^{1/2}$ denote the population inner product and 
the associated $L_2$ norm. 

%
%
\section{Multi-resolution group Lasso} 

In two subsections, we describe the multi-resolution structure for the approximation of the regression function in AM and introduce the MR-GL as the group Lasso estimator in the structure. 

\subsection{Multi-resolution structure} 
In a multi-resolution analysis, a nonparametric function $f_j(\cdot)$ can be 
written as a linear combination of basis functions in a block structure, 
\bel{f_j}
f_j(x) = \sum_{k=k_*}^{\infty} f_{j,k}(x),\quad f_{j,k}(x) 
= \sum_{\ell =1}^{2^{(k-1)\vee k_*}} \beta_{j,k,\ell}u_{j,k,\ell}(x), 
\eel
where 
the basis functions $u_{j,k,\ell}(\cdot), 1\le \ell\le 2^{(k-1)\vee k_*}, k\ge k_*$, are linearly independent 
in some infinite dimensional space for each $j$. 
We shall choose $\{u_{j,k,\ell}\}$ to form a complete system in the function spaces of interest to us, 
so that $f_j(x)$ can be approximated by finitely many terms in the series in \eqref{f_j} 
to an arbitrary level of accuracy. 
For simplicity and definitiveness, we take fixed $k_*\ge 0$ and block size $2^{(k-1)\vee k_*}$ 
in all infinite series expansions of the component functions $f_j(\cdot)$ throughout the paper. 
Under the observed design matrix $\bX$ with elements $x_{i,j}$, we write 
$\bff_j = \sum_{k=k_*}^{\infty} \bff_{j,k}$ and $\bff_{j,k}
= \sum_{\ell =1}^{2^{(k-1)\vee k_*}} \beta_{j,k,\ell}\bu_{j,k,\ell}$ 
with $\bu_{j,k,\ell}=(u_{j,k,\ell}(x_{1,j}),\ldots, u_{j,k,\ell}(x_{n,j}))^\top$. 

With the above block structure, we also consider truncated series of the form 
\bel{fbar_j}
\fbar_j(x) = \sum_{k=k_*}^{k^*} \fbar_{j,k}(x),\ 
\fbar_{j,k}(x) = \sum_{\ell=1}^{2^{(k-1)\vee k_*}} \betabar_{j,k,\ell}u_{j,k,\ell}(x),\ 
\betabar_{j,k,\ell}\in \R. 
\eel
We write $\bffbar_j = \sum_{k=k_*}^{\infty} \bffbar_{j,k}$ and 
$\bffbar_{j,k} = \sum_{\ell =1}^{2^{(k-1)\vee k_*}} \betabar_{j,k,\ell}\bu_{j,k,\ell}$ 
with the observed $\bX$.

Here and in the sequel, $f_j$ refers to a general candidate of the $j$-th AM component 
under consideration with functional version $f_j(\cdot)$ and realized vector version $\bff_j$, 
and $\fbar_{j}$ refers to a general candidate for its approximation up to a fixed 
resolution level $k^*$ with functional version $\fbar_{j}(\cdot)$ and vector version $\bffbar_{j}$. 
Similarly, $f_{j,k}$, $f_{j,k}(\cdot)$, $\bff_{j,k}$ and $\fbar_{j,k}$, $\fbar_{j,k}(\cdot)$, $\bffbar_{j,k}$ are the corresponding expressions at the resolution level $k$. 
Thus, 
\bel{F_{n,j}}
\scrF_{n,j}^{({\rm\tiny NP})} 
= \Big\{\bffbar_j: \fbar_j(\cdot) \hbox{ in } \eqref{fbar_j}\Big\} 
= \bigg\{\sum_{k=k_*}^{k^*}\sum_{\ell=1}^{2^{(k-1)\vee k_*}} b_{k,\ell}\bu_{j,k,\ell}: 
b_{k,\ell}\in\R\bigg\}
\eel
is the linear subspace of $\R^n$ generated by the approximation candidates in \eqref{fbar_j}. 
In this approximation scheme, the terms with $k>k^*$ are viewed as of ultra-high resolution 
and presumably ignorable. 

\def\Sobolev{{\rm\tiny Sobolev}}
The multi-resolution expansion \eqref{f_j} and the associated approximation \eqref{fbar_j} 
are essential in adaptive estimation of smooth functions in nonparametric analysis. 
For $\alpha >0$ define Sobolev-type norms 
$\|\bff_j\|_{\alpha,2,n}$ and $\|\bff_j\|_{\Sobolev,\alpha,n}$ by 
\bel{Sobolev-norm}
\|\bff_j\|_{\alpha,2,n}^2 = \sum_{k=k_*+1}^\infty 2^{2\alpha k}\|\bff_{j,k}\|_{2,n}^2,\   
\|\bff_j\|_{\Sobolev,\alpha,n}^2 = \|\bff_{j,k_*}\|_{2,n}^2+\|\bff_j\|_{\alpha,2,n}^2 
\eel
with the vector versions of $f_j$ and $f_{j,k}$ in \eqref{f_j}. 
We shall say that $f_j(\cdot)$ is $\alpha$-smooth (in the $\ell_2$ and empirical senses) 
when $\|\bff_j\|_{\alpha,2,n}$ is bounded 
and $f_j(\cdot)$ can be expresses in the infinite series expansion in \eqref{f_j}. 
Essentially, \eqref{Sobolev-norm} implies 
\bel{approx-error}
\big\| \bff_j - {\overline\bff}_j\big\|_{2,n} 
\le \|\bff_j\|_{\alpha,2,n}\bigg(\sum_{k>k^*}2^{-2\alpha k}\bigg)^{1/2}
= \frac{ 2^{-\alpha k^*}\|\bff_j\|_{\alpha,2,n}}{(4^\alpha -1)^{1/2}} 
\eel
for some $\bffbar_j\in \scrF_{n,j}^{({\rm\tiny NP})}$, for example with 
$\betabar_{j,k,\ell}=\beta_{j,k,\ell}$ in \eqref{fbar_j} but not necessarily. 
It is reasonable to assign the smoothness index $\alpha$ in \eqref{Sobolev-norm} as 
\bes
\E \|\bff_j\|_{\Sobolev,\alpha,n}^2 = O(1)\int\big\{f_j^2(x) + \big((d/dx)^\alpha f_j(x)\big)^2\big\}dx
\ees
when $x_{i,j}$ are random variables supported in a fixed compact interval 
with uniformly bounded marginal probability density functions 
and $u_{j,k,\ell}(\cdot)$ form the Fourier, wavelet or spline bases for each $j$. 
We note that $\E \|\bff_j\|_{\alpha,2,n}^2$ does not 
depend on $n$ when $x_{i,j}$ are identically distributed for the given $j$. 

Consider the case where the AM holds with only one non-zero component 
$f_j=f^*_j$ and iid noise $\veps_i\sim N(0,\sigma^2)$. 
Suppose that we are prepossessed to estimate $\bff^*=\E[\by|\bX]$ by the least squares 
projection of $\by$ to $\scrF_{n,j}^{({\rm\tiny NP})}$ when the oracular knowledge of 
$\bff^*=\bff_j^*$ is available. This oracle LSE has the average variance 
$2^{k^*}\sigma^2/n$ and maximum average squared bias of the order $2^{-2\alpha k^*}$ 
in view of \eqref{approx-error}. 
Thus, the $\ell_2$-regularization of $\bff_j$ as a group is not expected to be rate minimax unless 
$2^{k^*}\asymp n^{1/(2\alpha+1)}$, as it is not expected to outperform the oracle LSE 
in view of the minimax rate for fixed $p$ \citep{stone1985additive}. 
This forces the choice of $k^*$ to depend on $\alpha$ for rate optimality \citep{huang2010variable} 
and explains the general non-adaptivity in $\alpha$ in regularization schemes based on a fixed number of inner-product norms of $f_j(\cdot)$. 

The situation is different if we regularize each $f_j$ at individual resolution levels $k$, 
or equivalently regularize the effect represented by the $f_{j,k}$  in \eqref{f_j}. 
When a high-frequency $f^*_{j,k}$ carries little signal, 
its noisy expression in the data would be automatically thresholded or shrunk to reduce the 
variability of the estimator when $k_*\le k\le k^*$. 
However, when $k>k^*$, a great portion of the signal in the missing term 
$f_{j,k}$ is not recoverable even when the signal is large 
as it is excluded in the general approximation scheme, 
possibly causing excess bias and sub-optimality. 
These heuristics suggest the choice of a large $k^*$ to facilitate 
adaptation to the smoothness of the true $f^*_j(\cdot)$. 
The choice of the block size $2^{(k-1)\vee k_*}$ for the nonparametric $f_j$ 
is critical to the success of the MR-GL. 
First, it guarantees that different resolution levels are differentiated. 
Second, it groups large number of basis functions together to take advantage of the sparsity at 
the component function level. 

In the nonparametric regression setting, 
we still need the dimension of the approximation space $\scrF_{n,j}^{({\rm\tiny NP})}$ 
to be smaller than $n$ to relate it to its population version. 
As the dimension of the space 
is $2^{k^*}$ in \eqref{F_{n,j}},  
we take $2^{k^*} < n$ 
to capture most of the signal in 
$f^*_j$ when $f^*_j(\cdot)$ is $\alpha_j$-smooth with 
$n^{1/(2\alpha_j +1)}\le 2^{k^*}<n$, allowing all $\alpha_j>1/2$. 
Similarly, we shall not pick too large a $k_*$ 
as it forces the grouping of $2^{k_*}$ low frequency basis functions. 
For smooth functions, the baseline resolution level typically contains 
a fixed number of (constant or polynomial) terms. 
Later in our theoretical study, we will use $k_*\approx 2\log p$ to minimize the rate of 
theoretical error bounds. 
For notational convenience, we label the baseline resolution level $k_*$ 
so that there are totally $2^k$ basis functions up to level $k$. 

The AM also accommodate linear and group effects $f_j$ 
by allowing the first resolution in \eqref{f_j} and \eqref{fbar_j} 
to have flexible number of basis functions, 
\bel{ell^*_j}
f_{j,k_*}(x) = \sum_{\ell=1}^{d^*_j} \beta_{j,k_*,\ell}u_{j,k_*,\ell}(x),\quad 
\fbar_{j,k_*}(x) = \sum_{\ell=1}^{d^*_j} \betabar_{j,k_*,\ell}u_{j,k_*,\ell}(x), 
\eel
with $f_j=f_{j,k_*}$ and $\fbar_j=\fbar_{j,k_*}$ 
for a positive integer $d^*_j$, and $f_{j,k}(\cdot) = \fbar_{j,k}(\cdot) =0$ for all $k>k_*$.  
For example, when $f_j(x) = \beta_j x$ is a linear function, 
we take $\beta_{j,k_*,1}=\beta_j$ and $u_{j,k_*,1}(x)=x$ with $d^*_j=1$. 
Similarly, a group effect of dimension $d^*_j$ can be written as \eqref{ell^*_j} 
when $u_{j,k_*,\ell}(x_{i,j})$ are viewed as covariates in the group, $\ell\le d^*_j$. 
We shall call parametric components such linear and group effects in the AM. 

Putting the nonparametric components \eqref{f_j} and 
the parametric components \eqref{ell^*_j} together, the AM is written as 
\bel{K}
f = f(x_1,\ldots,x_p) = \sum_{j\in J_0} f_j(x_j) + \sum_{j\in J_1}\sum_{k=k_*}^\infty f_{j,k}(x_j)
= \sum_{(j,k)\in\scrK} f_{j,k}(x_j),
\eel
where $J_0$ and $J_1$ are respectively the index sets of the parametric and nonparametric 
components $f_j(\cdot)$, and 
$\scrK = \big\{(j,k): k=k_* \forall j\in J_0, k\ge k_* \forall j\in J_1\big\}$ 
is the index set for the parametric components represented at the nominal resolution level 
$k_*$ ($f_j=f_{j,k_*}$) and the nonparametric 
components represented at individual resolution levels $k\ge k_*$. 
In this multi-resolution structure, $f$ is approximated by 
\bel{K^*}
\fbar = \fbar(x_1,\ldots,x_p) = \sum_{j\in J_0} \fbar_j(x_j) + \sum_{j\in J_1}\sum_{k=k_*}^{k^*}
\fbar_{j,k}(x_j)
= \sum_{(j,k)\in\scrK^*} \fbar_{j,k}(x_j)
\eel
with $\fbar_{j,k}$ in \eqref{fbar_j} and \eqref{ell^*_j} and 
$\scrK^* = \big\{(j,k): k=k_* \forall j\in J_0, k_*\le k\le k^* \forall j\in J_1\big\}$. 
The noisy expression of the true $f^*_{j,k}$ in the data is still subject to regularization 
in the MR-GL for $k\le k^*$, but the MR-GL does not attempt to recover signals of nonparametric components in resolution levels $k > k^*$. 
We note that $k^*-k_*\asymp \log n$ typically in this multi-resolution scheme to allow adaptation 
to different levels of smoothness of the nonparametric components $f_j$. 

\subsection{Estimation through group Lasso}
In the vector notation, the AM \eqref{LM-0} can be written as 
\begin{equation}
    \by = \sum_{j=1}^p \bff_j + \bep,
    \label{LM-1}
\end{equation}
where $\by=(y_1,\ldots,y_n)^\top$ is the response vector, 
$\bff_j=(f_j(x_{1,j}),\ldots,f_j(x_{n,j}))^\top$ represents the contribution of the $j$-th covariate to the model, and $\bep=(\veps_1,\ldots,\veps_n)^\top$ is the noise vector.  

As the ultra high-frequency terms with $k>k^*$ are automatically excluded in 
the approximation scheme \eqref{K^*}, 
the MR-GL estimator of $\bff^*$ 
is defined as 
\bel{MR-GL-pred}
\hbf = \argmin_{\bff}\bigg\{\frac{\|\by - \bff\|_{2,n}}{2} + 
A_0\sum_{(j,k)\in\scrK^*} \lam_{j,k} \|\bff_{j,k}\|_{2,n}: \bff =\sum_{(j,k)\in\scrK^*} \bff_{j,k} \bigg\}, 
\eel
where the minimum is taken over the sum $\bff$ with $\bff_{j,k}=(f_{j,k}(x_{1,j}),\ldots,f_{j,k}(x_{n,j}))^\top$ 
corresponding to 
the $f_{j,k}(\cdot)$ in \eqref{f_j} for $j\in J_1$ and in \eqref{ell^*_j} for $j\in J_0$, 
$\lam_{j,k}$ are proper threshold levels, and $A_0 > 1$ is a constant, e.g. $A_0=1.01$ or $A_0=2$. 
As $\hbf$ is a regularized projection of $\by$ to a linear space 
through convex minimization in \eqref{MR-GL-pred}, it is uniquely defined. 

We may also use the basis vectors $\bu_{j,k,\ell}$ to write the MR-GL more explicitly. Let  
\bel{U_{j,k}}
\bff_{j,k} = \bU_{j,k}\bbeta_{j,k},\quad \bU_{j,k}\in\R^{d\times d_{j,k}},\quad \bbeta_{j,k}\in\R^{d_{j,k}}, 
\eel
where $\bU_{j,k}=(\bu_{j,k,1},\ldots,\bu_{j,k,d_{j,k}})\in \R^{n\times d_{j,k}}$,  
$d_{j,k}=d^*_j\,\forall j\in J_0$ as in \eqref{ell^*_j}, 
$d_{j,k}=2^{(k-1)\vee k_*}\,\forall j\in J_1$ as in \eqref{f_j} 
and $\bu_{j,k,\ell}=(u_{j,k,\ell}(x_{1,j}),\ldots, u_{j,k,\ell}(x_{n,j}))^\top$ . 
Recall that for the nonparametric components with $j\in J_1$, 
the first $\bU_{j,k}$, with $k=k_*$, represents the design matrix for a group of low-resolution 
components of $\bff_j$, while $\bU_{j,k}$ with $k>k_*$ blends in high-resolution components. 
Let 
\bel{d^*}
\bbetabar = (\bbetabar_{j,k}^\top,(j,k)\in\scrK^*)^\top\in \R^{d^*}
\eel 
be the coefficient vector of the $\fbar$ in \eqref{K^*} 
in the basis $\{u_{j,k,\ell}\}$ with $d^*  = \sum_{(j,k)\in\scrK^*}d_{j,k}$ 
and the $\scrK^*$ 
in \eqref{K^*}. 
The MR-GL estimates the coefficient vector \eqref{d^*} for an ideal approximate of 
$\bff^*$ by 
\bel{MR-GL-coef}
\hbbeta = \argmin_{\bb}
\bigg\{\frac{1}{2}\bigg\|\by - \sum_{(j,k)\in\scrK^*}\bU_{j,k}\bb_{j,k}\bigg\|_{2,n} + 
A_0\sum_{(j,k)\in\scrK^*} \lam_{j,k} \|\bU_{j,k}\bb_{j,k}\|_{2,n}\bigg\} 
\eel
with $\hbbeta= (\hbbeta^\top_{j,k},(j,k)\in\scrK^*)^\top$ and $\bb = (\bb^\top_{j,k},(j,k)\in\scrK^*)^\top$. This can be viewed as group Lasso in linear regression 
when the columns $\bu_{j,k,\ell}$ of $\bU_{j,k}$ are treated as design vectors 
and $\bff_{j,k}$ are treated as group effects in \eqref{U_{j,k}}. 
While \eqref{MR-GL-coef} is not necessarily unique, any solution of it 
would recover the unique solution of \eqref{MR-GL-pred} via 
\bel{MR-GL-coef-pred}
\hbf = \sum_{j=1}^p \hbf_j = \sum_{(j,k)\in\scrK^*}\hbf_{j,k},\quad 
\hbf_j = \sum_{k: (j,k)\in\scrK^*}\hbf_{j,k},\quad 
\hbf_{j,k} = \bU_{j,k}\hbbeta_{j,k}, 
\eel
with the $\bU_{j,k}$ in \eqref{U_{j,k}}.  
The vectors $\hbf_j$ and $\hbf_{j,k}$ have the interpretation 
as the estimated AM components and their representation in individual resolution levels. 
Thus, (\ref{MR-GL-pred}), \eqref{MR-GL-coef} and \eqref{MR-GL-coef-pred} are equivalent and 
can be all viewed as MR-GL. 
The MR-GL enjoys the advantage of group Lasso, a computationally feasible and parsimonious solution, 
and adaptation to the smoothness of the component functions through group selection 
in the multi-resolution decomposition. 

Finally we discuss the choice of the thresholding level $\lam_{j,k}$. 
As \eqref{MR-GL-coef} is group Lasso, the literature provides a variety of its error bounds 
for $\hbf - \bff^*$ when 
\bel{lam-gen}
\P\big\{\Omega_0\big\} = 1+o(1)\ \hbox{ with }\ 
\Omega_0 = \bigg\{\max_{(j,k)\in\scrK^*}\frac{\|\bP_{j,k}(\by - \bff^*)\|_2}{\lam_{j,k}} \le 1\bigg\},  
\eel
where $\bP_{j,k}$ is the orthogonal projection from $\R^n$ to the 
column space of $\bU_{j,k}$ in \eqref{U_{j,k}}.  
When $\by - \bff^*$ has iid $N(0,\sigma^2)$ components, \eqref{lam-gen} holds when 
\bel{lam}
\lam_{j,k} 
= \begin{cases} 
\sigma\big((d^*_j/n)^{1/2}+\sqrt{(2/n)\log(p/\eps)}\big), & j\in J_0,  \cr 
\sigma\big((2^k/n)^{1/2}+\sqrt{(2/n)\log(p/\eps)}\big), & j\in J_1, 
\end{cases} 
\eel
with $\eps\in (0,1]$. 
The above $\lam_{j,k}$ also provides \eqref{lam-gen} when $\by - \bff^*$ is a sub-Gaussian vector 
with a slightly inflated $\sigma$ \citep{huang2010benefit}. 
Throughout the paper, we use $f^*_j$ and $f^*_{j,k}$ to denote the estimation target which could be the true $f_j$ and $f_{j,k}$ associated with $\bff^*=\E[\by|\bX]$ or any $\bff^*$ satisfying \eqref{lam-gen} via 
$\bff^*=\sum_{j=1}^p \bff^*_j = \sum_{(j,k)\in\scrK^*}\bff^*_{j,k}$. 

\section{Theoretical Results for Fixed Designs}\label{3-1} 
In this section, we carry out a theoretical analysis of the MR-GL in both fixed  
and random design settings. 
The main results are summarized by oracle inequalities stated 
in Theorem \ref{th-1} through Theorem \ref{th-5}. 
In the fixed design setting, we prove that the MR-GL adaptively achieves and generalizes 
existing rate optimal error bounds for various classes of the unknown $\bff^*$
under an empirical compatibility condition. 
In the random design setting, we provide oracle inequalities in both the empirical and population 
error measures as well as sufficient condition on the sample size 
to guarantee the empirical compatibility condition under mild conditions on the design population. 

In this section, we consider fixed designs in which 
the covariates $x_{i,j}$ in \eqref{LM-0} and \eqref{LM-1} are treated as deterministic. 
We measure the performance of the MR-GL by the in-sample squared prediction error  
$\|\hbf-\bff^*\|_{2,n}^2$ where $\bff^*$ is the estimation target. 
Typically $\bff^*=\E[\by]$ as $\bX$ is deterministic but 
our theorems also apply to any $\bff^*$ satisfying \eqref{lam-gen}. 
As the MR-GL \eqref{MR-GL-pred} can be viewed as group Lasso via \eqref{MR-GL-coef} and 
\eqref{MR-GL-coef-pred}, our results are closely related to 
the existing group Lasso theory  
\citep{nardi2008asymptotic, huang2010benefit, lounici2011oracle, buhlmann2011statistics, negahban2012unified, mitra2016benefit}. 
However, due to the choice of $\lam_{j,k}$ in \eqref{lam} and 
the truncation of the nonparametric components  
beyond the resolution level $k^*$, we need different error bounds 
based on the approximate sparsity of the group components. 
We discuss in separate subsections oracle inequalities for the MR-GL in general, 
the adaptive minimaxity of the MR-GL to the smoothness index $\alpha$ in \eqref{Sobolev-norm} 
in the nonparametric AM, 
and the connection of our results to the group Lasso in linear regression. 

\subsection{General error bounds}\label{3-1-1} 
For the $\bU_{j,k}$ in \eqref{U_{j,k}}, 
$\bb =(\bb_{j,k}^\top,(j,k)\in\scrK^*)^\top\in\R^{d^*}$ as in \eqref{MR-GL-coef} 
and $S\subset \scrK^*$, define 
\bel{pen_S}
\pen_S(\bb)=\sum_{(i,k)\in S}\lam_{j,k}\|\bU_{j,k}\bb_{j,k}\|_{2,n}
\eel
as the partial penalty on $\bb$ in the index set $S$. 
For $\xi\ge 1$ and $S\subset \scrK^*$, define 
\bel{C-pred-1}
C_{\rm pred}(\xi,S) 
= \sup_{\bb} \frac{\big\{\pen_S(\bb) - \xi^{-1}\pen_{S^c}(\bb)\big\}_+^2}
{\|\lam_S\|_2^2\|\sum_{(j,k)\in\scrK^*}\bU_{j,k}\bb_{j,k}\|_{2,n}^2}
\eel
as a constant factor for empirical prediction error bounds, 
with $S^c = \scrK^*\setminus S$, $\|\lam_S\|_2=(\sum_{(j,k)\in S}\lam_{j,k}^2)^{1/2}$ 
and the convention $0/0 = 0$. 
This prediction factor 
is no greater than the reciprocal of the compatibility coefficient 
for group Lasso \citep{buhlmann2011statistics} 
as we will discuss in Section \ref{3-1-3}.

\begin{theorem}\label{th-1} 
{\rm (i)} Let $\hbf$ be as in \eqref{MR-GL-coef-pred} with 
$\hbbeta=(\hbbeta_{j,k}^\top,(j,k)\in\scrK^*)^\top\in\R^{d^*}$ in \eqref{MR-GL-coef}. 
Let $\bffbar = \sum_{(j,k)\in\scrK^*} \bffbar_{j,k}$ with $ \bffbar_{j,k} = \bU_{j,k}\bbetabar_{j,k}$,   
$\bU_{j,k}\in \R^{n\times d_{j,k}}$ in \eqref{U_{j,k}} and 
$\bbetabar=(\bbetabar_{j,k}^\top,(j,k)\in\scrK^*)^\top$ in \eqref{d^*}. 
Let $\bP_{j,k}$ be the orthogonal projection to the range of $\bU_{j,k}$, 
$\pen_S(\bb)$ in \eqref{pen_S} and $C_{\rm pred}(\xi,S)$ in \eqref{C-pred-1} 
with $\xi = (A_0+1)/(A_0-1)>1$ and $S\subset \scrK^*$. 
Then, 
\bel{th-1-1}
\big\|\hbf - \bff^*\big\|_{2,n}^2\le B_S+\Deltabar_S,\quad 
\big\|\hbf - \bffbar\big\|_{2,n}^2 + \big\|\hbf - \bff^*\big\|_{2,n}^2\le 4B_S+2\Deltabar_S, 
\eel
with $\Deltabar_S= \big\|\bffbar- \bff^*\big\|_{2,n}^2 + 4A_0\pen_{S^c}(\bbetabar)$ and 
$B_S = (A_0+1)^2C_{\rm pred}(\xi,S)\|\lam_S\|_2^2$, when 
\bel{th-1-2}
\|\bP_{j,k}(\by - \bff^*)\|_{2,n} \le \lam_{j,k}\quad \forall\ (j,k)\in \scrK^*. 
\eel
In particular, for $S=\{(j,k)\in\scrK^*: \|\bffbar_{j,k}\|_{2,n}\ge A_0\lam_{j,k}\}$, 
\bel{th-1-3}
&& \big\|\hbf - \bffbar\big\|_{2,n}^2 + \big\|\hbf - \bff^*\big\|_{2,n}^2 
\cr &\le& 2\big\|\bffbar- \bff^*\big\|_{2,n}^2+C^*_{\rm pred}(\xi,S)\textsum_{(j,k)\in\scrK^*}
\lam_{j,k}^2\wedge\big(\lam_{j,k}\|\bffbar_{j,k}\|_{2,n}\big)
\eel
with $C^*_{\rm pred}(\xi,S)=\max\big\{8A_0^2, 4(A_0+1)^2C_{\rm pred}(\xi,S)\big\}$ 
and $\bffbar_{j,k}= \bU_{j,k}\bbetabar_{j,k}$. \smallskip

\noindent
{\rm (ii)} Suppose $\by - \bff^*$ has iid $N(0,\sigma^2)$ entries. Then, \eqref{th-1-2} 
holds with probability at least $1-\eps/\sqrt{2\log(p/\eps)}$ when $\lam_{j,k}$ are given by 
\eqref{lam} with $0 < \eps \le 1\le 2\log(p/\eps)$. 
\end{theorem}

In the above theorem, $\bff^*$ can be viewed as $\E[\by]$ as the design is treated as deterministic. 
However, \eqref{th-1-1} and \eqref{th-1-3} apply to any vector $\bff^*\in \R^n$ satisfying \eqref{th-1-2} 
with high probability. In particular, $\bff^*$ is not required to have additive components. 
In Theorem \ref{th-1}, \eqref{th-1-3} bounds the in-sample squared prediction error  
$\|\hbf-\bff^*\|_{2,n}^2$ of the MR-GL by the sum of the squared approximation 
error $\|\bffbar- \bff^*\|_{2,n}^2$ of a deterministic approximation candidate $\bffbar$ 
generated by a function $\fbar$ in \eqref{K^*} 
and the product of the prediction factor $C^*_{\rm pred}(\xi,S)$ and a normalized 
complexity measure of $\bffbar$. 
Thus, Theorem \ref{th-1} asserts that the in-sample prediction error of the MR-GL is small 
when the estimation target $\bff^*$ can be approximated by some unknown sparse $\bffbar$, 
and the sparsity of $\bffbar$ means the sparsity of its functional components as well as the decay of the signal strength at high resolution levels due to the smoothness of the individual nonparametric components, respectively in the index $j$ and in the index $k$ given $j$. 

Here $\bffbar$ can be any member of the approximate space 
\bel{F_n}
\scrF_n &=& \bigg\{\bffbar = \sum_{j=1}^p \bffbar_j: 
\bffbar_j = \bffbar_{j,k_*}\ \forall\, j\in J_0,\ 
\bffbar_j\in \scrF_{n,j}^{({\rm\tiny NP})}\ 
\forall\, j\in J_1\bigg\}
\cr &=& \bigg\{\bffbar = \bU\bbetabar = \sum_{(j,k,\ell)\in\scrL^*} \betabar_{j,k,\ell}\bu_{j,k,\ell}: 
\bbetabar\in\R^{d^*}, \betabar_{j,k,\ell}\in\R\bigg\}, 
\eel
where $\scrF_{n,j}^{({\rm\tiny NP})}$ is the approximation space in \eqref{F_{n,j}} for 
the nonparametric components, 
$\scrL^*=\{(j,k,\ell): (j,k)\in \scrK^*, 
\ell\le d^*_j\,\forall j\in J_0, \ell\le 2^{(k-1)\vee k_*}\,\forall j\in J_1\}$, $d^*=|\scrL^*|$, 
$\bU$ is the $n\times d^*$ matrix composed of columns 
$\bu_{j,k,\ell}=(u_{j,k,\ell}(x_{1,j}),\ldots, u_{j,k,\ell}(x_{n,j}))^\top$ with $(j,k,\ell)\in\scrL^*$, 
and $\bbetabar$ is the coefficient vector as in \eqref{d^*}. 
We note that when $d^*=|\scrL^*| > n$, $\scrF_n$ typically fills the entire space $\R^n$ 
but the approximation candidates $\bffbar$ is still meaningful through its representation 
$\bU\bbetabar$ as in sparse linear regression. Moreover, the coefficients allow 
identification of the components $\bffbar_j$ through the functional representations of 
$\fbar_j$ in \eqref{fbar_j} and \eqref{ell^*_j}. 

The second term on the right-hand side of \eqref{th-1-3} can be viewed as a normalized 
complexity measure of the approximation candidate $\bffbar$ in the following sense. 
Let $\delta_k=I\{k>k_*\}$ and $d_{j,k}$ be as in \eqref{U_{j,k}}.
For each group with index $(j,k)$ in $\scrK^*$, the quantity 
\bel{complexity}
\frac{\lam_{j,k}^2\wedge\big(\lam_{j,k}\|\bffbar_{j,k}\|_{2,n}\big)}{\sigma^2/n} 
= \big(2^{\delta_k/2}d_{j,k}^{1/2}+\sqrt{2\log(p/\eps)}\big)^2
\min\bigg(1,\frac{\|\bffbar_{j,k}\|_{2,n}}{\lam_{j,k}}\bigg). 
\eel
is roughly the number of data points needed to estimate $\bffbar_{j,k} = \bU_{j,k}\bbetabar_{j,k}$ 
as expressed in \eqref{U_{j,k}}, as $d_{j,k}$ is the number of columns of $\bU_{j,k}$.  
For strong signals $\|\bffbar_{j,k}\|_{2,n}\ge \lam_{j,k}$, 
the nominal degrees of freedom $d_{j,k}$ is inflated to \eqref{complexity} 
to take into account of the uncertainty about its signal strength in group selection. 
For weaker signals, \eqref{complexity} discounts the complexity by the ratio between 
the signal strength $\|\bffbar_{j,k}\|_{2,n}$ and penalty level $\lam_{j,k}$. 

The above connection between the error bound \eqref{th-1-3} and the complexity of 
the approximation of the unknown 
indicates its rate optimality in a broad range of settings, and thus the adaptive 
optimality of the MR-GL as we will discuss below. 

\subsection{Adaptive optimality in nonparametric AM}\label{3-1-2} 
The main objective of this subsection is to 
present the implications of Theorem \ref{th-1} in the nonparametric AM in which all 
the components $f_j$ are nonparametric, i.e. $J_0=\emptyset$ and $J_1=\{1,\ldots,p\}$, 
while the discussions below up to the statement of Corollary~\ref{cor-1}, 
including Theorem \ref{th-2}, are applicable 
to the general semi-parametric setting with possibly nonempty $J_0$.  
We shall focus on the theory in which the complexity of $f_j$ is measured by the Sobolev-type 
norm \eqref{Sobolev-norm} with a common index $\alpha$ and possibly a second 
smaller common index $\alpha_0$ to facilitate more direct comparisons between our results and those in the existing literature. 

We first connect the second term
in the error bound \eqref{th-1-3} to the Sobolev-type norm  
$\|\bff_j\|_{\alpha,2,n}$ in \eqref{Sobolev-norm}. 
For real $c$ and $0\le q\le 1$, define 
\bes
J_{c}^{(q)}(k_1,k_2) 
&=& \begin{cases}\Big(\sum_{k=k_1+1}^{k_2} 2^{c k/(1-q/2)}\Big)^{1-q/2}, & c \le 0, 
\cr \Big(\sum_{k=0}^{k_2-k_1-1} 2^{-c k/(1-q/2)}\Big)^{1-q/2}, & c > 0,
\end{cases}
\ees
for all integers $0\le k_1\le k_2$, with $J_{c}^{(q)}(k,k)=0$. 
Consider two smoothness indices $0\le\alpha_0\le\alpha$. 
By the definition of the norm $\|\bffbar_j\|_{\alpha,2,n}$ in \eqref{Sobolev-norm}, 
the sum is bounded by 
\bes
\sum_{k=k_*+1}^{k^*}\lam_{j,k}^2\wedge\big(\lam_{j,k}\|\bffbar_{j,k}\|_{2,n}\big)
\ge \max_{k_*<k\le k^*} \lam_{j,k}\min\bigg(\lam_{j,k},
\frac{\|\bffbar_j\|_{\alpha_0,2,n}}{2^{\alpha_0 k}},\frac{\|\bffbar_j\|_{\alpha,2,n}}{2^{\alpha k}}\bigg). 
\ees
When $\lam_{j,k}\asymp 2^{k/2}$ with $2^{k_*}\asymp \log(p/\eps)$ in \eqref{lam}, 
the above lower bound is essentially sharp as the summations below and above the maximizing  
$k$ are both bounded by geometric series of the form $J_{c}^{(0)}(k_1,k_2)$. 
The following proposition refines the above argument and facilitates the use of 
the sequence norms of $\|\bffbar_j\|_{\alpha,2,n}$ and $\|\bffbar_j\|_{\alpha_0,2,n}$ to bound 
the sum $\sum_{j\in J_1}\sum_{k=k_*+1}^{k^*}\lam_{j,k}^2\wedge\big(\lam_{j,k}\|\bffbar_{j,k}\|_{2,n}\big)$ 
in \eqref{th-1-3}. 

\begin{proposition}\label{prop-1}
Let $\alpha \ge 1/2$, $0\le q\le 1$ and $0\le\alpha_0\le\alpha$. 
Define
\bel{gamma}
\gamma  = \frac{(2-q)(\alpha-1/2) + (1 - q/2 - q\alpha_0)_+}{\alpha - 1/2 + (1 - q/2 - q\alpha_0)_+}
\eel
with $\gamma=1\,$ for $\alpha = 1/2$. 
Let $\sigma_n>0$ and $\lam_k \le \sigma_n 2^{k/2}$ for $k_*<k\le k^*$. Then, 
\bel{prop-12} 
&& \sum_{k=k_*+1}^{k^*}\lam_{k}\min\Big(\|\bffbar_{j,k}\|_{2,n},\lam_{k}\Big)
\le \sigma_n^{\gamma}J_{q,\alpha,\alpha_0}(k_*,k^*)\big\{\|\bffbar_{j}\|_{\alpha_0,2,n}^q\big\}^{1-\rho}
\|\bffbar_{j}\|_{\alpha,2,n}^{\rho},\quad
\eel
where $\rho = (1-q/2-\alpha_0q)_+/\{\alpha - 1/2 + (1-q/2-\alpha_0q)_+\}$ 
and 
\bes
J_{q,\alpha,\alpha_0}(k_*,k^*) 
=\big\{(1/\rho-1)^{\rho} + (1/\rho-1)^{\rho-1}\big\} 
\big\{J_{1-q/2-\alpha_0q}^{(q)}(k_*,k^*)\big\}^{1-\rho}
\big\{J_{\alpha-1/2}^{(1)}(k_*,k^*)\big\}^{\rho}. 
\ees 
\end{proposition} 

\begin{remark}\label{remark-1}
We note that $1\le \gamma \le \min\{4\alpha/(2\alpha+1), 2-q\}$ always holds, 
\bes
\gamma = \begin{cases}
\min\{4\alpha/(2\alpha+1), 2-q\}, & \alpha_0=\alpha, 
\cr 4\alpha/(2\alpha+1), & q=0, 
\cr (2-q)2\alpha/(2\alpha+1-q), & \alpha_0=0, 
\cr 1, & \alpha=1/2\ \hbox{ or }\ q=1, 
\end{cases}
\ees
$\gamma < 2-q$ iff $1-q/2-\alpha_0q>0$, 
and $\gamma < 4\alpha/(2\alpha+1)$ if $q\neq 0$ and $\alpha_0<\alpha$.  
Moreover, 
\bes
\rho = \begin{cases} (2-q-\gamma)/(1-q), & 0\le q<1, 
\cr 2/(2\alpha+1), & q=0, 
\cr (2-q)/(2\alpha+1-q), & \alpha_0=0, 
\cr 1, & \alpha=1/2, 
\cr (1-1/2-\alpha_0)_+/(\alpha-1/2+(1-1/2-\alpha_0)_+), & q=1, 
\end{cases}
\ees 
and $\gamma + q(1-\rho)+\rho=2$ for the $\rho$ in \eqref{prop-12}.  
\end{remark}

\begin{remark} 
We note that 
$J_{c}^{(q)}(k_1,k_2)\asymp \min\big\{|c|^{q/2-1},(k_2-k_1)^{1-q/2}\big\}$ uniformly in $(c,q,k_1,k_2)$. 
For $c \neq 0$, $J_{c}^{(q)}(k_1,k_2)\le \big\{1-2^{-|c|/(1-q/2)}\big\}^{q/2-1}$ in Theorem~\ref{th-1}. However, this upper bound is not accurate when $c$ is close to zero. 
\end{remark}

\medskip
We need to specify an $\bffbar$ in Theorem \ref{th-1} so that the first term on the right hand side of (\ref{th-1-3}) is small. 
When $\bff^*$ arises from the AM \eqref{K}, 
a natural way of constructing $\bffbar$ is to truncate the ultra-high-frequency 
terms of $\bff^*$, or equivalently to set the coefficients $\betabar_{j,k,\ell}$ in \eqref{fbar_j} 
equal to the true version. 
Formally, suppose the AM holds for the true regression function in the sense of 
$f_j(x)=f_j^*(x)$ and $\E[\veps_i]=0$ in \eqref{LM-0}. When $f^*_j(x)$ has an infinite series expansion 
of the form \eqref{f_j} with coefficients $\beta^*_{j,k,\ell}$, we may simply set 
$\betabar_{j,k,\ell} = \beta^*_{j,k,\ell}$ in \eqref{fbar_j}. This gives 
$\bffbar = 
\sum_{(j,k)\in\scrK^*}\bffbar_{j,k}\in \scrF_n$ as in \eqref{F_n} with 
\bel{bffbar}
\bffbar_{j,k} = \bff^*_{j,k}\ \forall\, (j,k)\in \scrK^*,\quad 
\bff_{j,k}^*=\sum_{\ell=1}^{d_{j,k}}\beta^*_{j,k,\ell}\bu_{j,k,\ell}\ \forall\, (j,k)\in \scrK,
\eel
where $\bu_{j,k,\ell}=(u_{j,k,\ell}(x_{1,j}),\ldots, u_{j,k,\ell}(x_{n,j}))^\top$, $d_{j,k}$ are as in \eqref{U_{j,k}}, 
and $\scrK$ and $\scrK^*$ 
are as in \eqref{K} and  \eqref{K^*} respectively. 
When the nonzero components of $\bff^*-\bffbar$ are not highly correlated, we expect 
that \eqref{bffbar} would provide 
\bel{approximate1}
\|\bff^* - \bffbar\|_{2,n}^2 \lesssim \sum_{j\in J_1}\bigg\|\sum_{k=k^*+1}^\infty  \bff^*_{j,k}\bigg\|_{2,n}^2
\le \frac{2^{-2\alpha k^*}}{4^\alpha-1}\sum_{j\in J_1}\|\bff^*_j\|_{\alpha,2,n}^2
\eel
where $\|\bff_j\|_{\alpha,2,n}$ is the Sobolev-type norm defined in (\ref{Sobolev-norm}) 
and $J_1$ is as in \eqref{K}. 
Thus, the first term on the right hand side of (\ref{th-1-3}) is controlled by 
the $\ell_2$ norm of $\|\bff_j^*\|_{\alpha,2,n}$. 
Alternatively we may use the following cruder bound which always holds: 
\bel{approximate2}
\|\bff^* - \bffbar\|_{2,n} \le \sum_{j\in J_1}\bigg\|\sum_{k=k^*+1}^\infty  \bff^*_{j,k}\bigg\|_{2,n} 
\le \frac{2^{- \alpha k^*}}{(4^\alpha-1)^{1/2}}\sum_{j\in J_1}\|\bff^*_j\|_{\alpha,2,n}. 
\eel
Thus, the first term on the right hand side of (\ref{th-1-3}) is explicitly controlled by 
the $\ell_1$ norm of $\|\bff_j^*\|_{\alpha,2,n}$ without condition on the design.  

We are now ready to present the adaptive optimality of the MR-GL. Define 
\bel{M_alpha_n}
M_{\alpha,q,n}^q=\sum_{j\in J_1}\|\bff^*_j\|_{\alpha,2,n}^q,\ 
M_{q,n}^{q,\BR}=\sum_{j=1}^p (\lam_{j,k_*}/\lam_0)^{2-q}\|\bffbar_{j,k_*}\|_{2,n}^q, 
\eel
with $\lam_0=\sigma\sqrt{(2/n)\log(p/\eps)}$,
respectively as the $q$-power of the $\ell_q$ "norm" of the norms 
$\|\bff_j^*\|_{\alpha,2,n} = \{\sum_{k=k_*+1}^\infty 2^{2\alpha k}\|\bff_{j,k}^*\|_{2,n}^2\}^{1/2}$ 
in \eqref{Sobolev-norm} for the nonparametric components 
and the $q$-power of a weighted $\ell_q$ norm of the 
$\ell_{2,n}$ norms of $\bff_{j,k_*}^*$ for all parametric and nonparametric 
components, for all $q \ge 0$, with $q=0$ treated as the limit at $q=0+$.
While $M_{q,n}^{q,\BR}$ describes the complexity of the representation of all the components 
$f_j$ at the baseline resolution level $k=k_*$, $M_{\alpha,q,n}^q$ measures the complexity 
of the nonparametric components $\{f_j, j\in J_1\}$ beyond the baseline resolution.  
In the nonparametric AM with $J_1=\{1,\ldots,p\}$, 
$\lam_{j,k_*}/\lam_0 = 2^{k_*/2}/\sqrt{2\log(p/\eps)}+1 \le 3$.

\begin{theorem}\label{th-2} 
Suppose $\by - \bff^*$ has iid $N(0,\sigma^2)$ entries with 
$\bff^*=\sum_{j=1}^p\sum_{k\ge k_*}\bff_{j,k}^*$ 
where $\bff^*_{j,k}$ are as in (\ref{bffbar}). 
Let $k^*\ge k_*$ be integers satisfying $2^{k_*-1} < 2\log(p/\eps) \le 2^{k_*}$ 
and $2^{k^*}\ge n^{1/(2\alpha_*+1)}$ for some $\alpha_*>0$. 
Let $\hbf$ be the MR-GL estimator in \eqref{MR-GL-coef-pred} with the estimated 
coefficients $\hbbeta_{j,k}$ in \eqref{MR-GL-coef}, 
$\lam_{j,k}= \sigma_n\big(2^{k/2}+\sqrt{2\log(p/\eps)}\big)$ with $\sigma_n=\sigma/n^{1/2}$,  
and constant $A_0>1$. 
Let $\bffbar$ be as in \eqref{bffbar}, 
$\{q,q_0\}\subset [0,1]$, $0\le \alpha_0\le \alpha$ with $\alpha>1/2$, 
$\gamma$ as in \eqref{gamma}, $\rho$ as in \eqref{prop-12} and 
$q(1-\rho)/q_2+\rho/q_1 = 1$ with $q_1\ge\rho$. 
Then, 
\bel{th-2-1}
&& \big\|\hbf - \bffbar\big\|_{2,n}^2 + \big\|\hbf - \bff^*\big\|_{2,n}^2 
\\ \nonumber &\le& \frac{n^{- 2\alpha/(2\alpha_*+1)}M_{\alpha,1,n}^2}{(4^{\alpha}-1)/2} 
+C^*_{\rm pred}(\xi,S)\bigg[\sigma_n^\gamma 
\frac{M_{\alpha_0,q_2,n}^{q(1-\rho)}M_{\alpha,q_1,n}^{\rho}}{2^{-2}J^{-1}_{q,\alpha,\alpha_0}(k_*,k^*)}
+ \lam_0^{2-q_0}M^{q_0,\BR}_{q_0,n}\bigg], 
\eel
with at least probability $1- \eps/\sqrt{2\log(p/\eps)}$ 
and with the option of replacing $M_{\alpha,1,n}^2$ by $O(M_{\alpha,2,n}^2)$ in the first term 
when \eqref{approximate1} holds, 
where 
$\lam_0=\sigma_n\sqrt{2\log(p/\eps)}$, 
$J_{q,\alpha,\alpha_0}(k_*,k^*)$ is as in \eqref{prop-12} and 
$C^*_{\rm pred}(\xi,S)$ as in \eqref{th-1-3}. 
In particular, the following hold when $C^*_{\rm pred}(\xi,S)=O(1)$, 
\eqref{approximate1} holds and $\alpha_0\neq 1/q-1/2$. \\
{\rm (i)} With $q=0$, $\gamma=4\alpha/(2\alpha+1)=2-\rho$, $q_1=\rho$ and 
$\#\{j: \bff^*_j\neq {\bf 0}\} \le s_0$ in \eqref{th-2-1}, 
\bel{th-2-2}
    \|\hbf - \bff^*\|_{2,n}^2 \lesssim  
    s_0n^{- 2\alpha/(2\alpha_*+1)}M_{\alpha,\infty,n}^2
    + s_0\sigma_n^{4\alpha/(2\alpha+1)}M_{\alpha,\infty,n}^{2/(2\alpha+1)} 
    + s_0\lam_0^2. 
\eel
{\rm (ii)} With $\alpha_0=\alpha$, $\gamma=\min\{4\alpha/(2\alpha+1), 2-q\}$, 
$q_1=1$ and $q_2=q$ in \eqref{th-2-1}, 
\bel{th-2-3}
    \|\hbf - \bff^*\|_{2,n}^2 \lesssim n^{- 2\alpha/(2\alpha_*+1)}M_{\alpha,2,n}^2
    + \sigma_n^{\gamma}M_{\alpha,q,n}^{2-\gamma}
    + \lam_0^{2-q_0}M_{q_0,n}^{q_0,\BR}. 
\eel
{\rm (iii)} With $\alpha_0=0$, $\gamma=(2-q)2\alpha/(2\alpha+1-q)$, 
$\rho=\gamma/(2\alpha)$, 
$q_1=1$ and $q_2=q$ in \eqref{th-2-1}, 
\bel{th-2-4}
    \|\hbf - \bff^*\|_{2,n}^2 \lesssim 
    n^{- 2\alpha/(2\alpha_*+1)} M_{\alpha,2,n}^2
    + \sigma_n^{\gamma}M_{0,q,n}^{q(1-\rho)}M_{\alpha,1,n}^{\rho}
    + \lam_0^{2-q_0}M_{q_0,n}^{q_0,\BR}. 
\eel
\end{theorem}


The oracle inequalities in the above theorem use three terms to bound 
the in-sample squared prediction error of the MR-GL. 
The first term represents the ultra high resolution part of $f^*_j$ ignored by the MR-GL. 
While $\alpha_*>0$ depends on the tuning parameter $k^*$ 
via $2^{k^*}\ge n^{1/(2\alpha_*+1)}$, 
we can typically take a small $\alpha_*\in (0,1/2)$ 
so that $\alpha >\alpha_*$ and it would be reasonable to expect 
\bel{th-2-1a} 
n^{- 2\alpha/(2\alpha_*+1)}M_{\alpha,1,n}^2
=o(1)\sigma_n^\gamma C^*_{\rm pred}(\xi,S)  J_{q,\alpha,\alpha_0}(k_*,k^*) 
M_{\alpha_0,q_2,n}^{q(1-\rho)}M_{\alpha,q_1,n}^{\rho}
\eel
in \eqref{th-2-1} and to remove the first terms in \eqref{th-2-2}, \eqref{th-2-3} and \eqref{th-2-3}.   
For example, for $\alpha=2$ and 
$\alpha_*\le 1/2$, $2\alpha/(2\alpha_*+1) \ge 2$ 
so that \eqref{th-2-2} holds without requiring \eqref{approximate1} 
when $M_{\alpha,\infty,n} = O(1)$ 
because $n^{- 2\alpha/(2\alpha_*+1)}M_{\alpha,1,n}^2 
\le s_0 n^{-2} 
M_{\alpha,\infty,n}^2 
\le n^{-1}M_{\alpha,\infty,n}^2$ 
when $s_0\le n$. 

The third term on the right-hand sides of the oracle inequalities in Theorem \ref{th-2} 
represents the risk of estimating the $f_j$ at the baseline resolution level $k=k_*$ 
with model uncertainty adjustment. 
When $d^*_j/\log(p/\eps)$ is uniformly bounded as in the nonparametric AM, this 
term is of the same order as the in-sample squared 
prediction error rate for the Lasso in linear regression. 
In general, this term matches the squared error rate for the group Lasso. 

As the first term is typically of smaller order and the third term is unavoidable, 
the optimality of the oracle inequalities in Theorem \ref{th-2} is 
largely determined by the second term on the right-hand side. 
When $\sigma_n^\gamma\asymp n^{-2\alpha/(2\alpha+1)}$, i.e. $\gamma = 4\alpha/(2\alpha+1)$, 
their rate optimality is evident because the rate matches the minimax risk for the estimation of a single smooth 
function. The rate optimality of the MR-GL 
in several special cases of the nonparametric AM 
are discussed below and compared with the existing results. 

\begin{corollary}\label{cor-1}  
Consider the nonparametric AM with $J_1=\{1,\ldots,p\}$, 
$\alpha\ge\alpha_*$ and $\log(1/\eps) = O(\log p)$, for example $\eps = p^{-a_0}$ with fixed $a_0\ge 0$. \\
{\rm (i)} With $M_{\alpha,\infty,n}=O(1)$ in \eqref{th-2-2},  
\bel{cor-1-1}
    \|\hbf - \bff^*\|_{2,n}^2 \lesssim  
    s_0n^{-2\alpha/(2\alpha+1)} + s_0(\log p)/n. 
\eel
{\rm (ii)} With $q_0=q$ and $M_{\alpha,q,n}^q\vee M_{q,n}^{q,\BR}=O(1)$ in \eqref{th-2-3}, 
\bel{cor-1-2}
    \|\hbf - \bff^*\|_{2,n}^2 \lesssim 
    n^{- 2\alpha/(2\alpha+1)}
    + \big((\log p)/n\big)^{1-q/2}. 
\eel
{\rm (iii)} With $q_0=q$ and $M_{\alpha,1,n}\vee M_{0,q,n}^q\vee M_{q,n}^{q,\BR}=O(1)$ 
in \eqref{th-2-4}, 
\bel{cor-1-3}
    \|\hbf - \bff^*\|_{2,n}^2 \lesssim 
    n^{- (2-q)\alpha/(2\alpha+1-q)}
    + \big((\log p)/n\big)^{1-q/2}. 
\eel
{\rm (iv)} With $q=q_0=0$ and $\#\{j: \bff^*_j\neq {\bf 0}\} \le s_0$ in \eqref{th-2-4}, 
\bel{cor-1-4}
    \|\hbf - \bff^*\|_{2,n}^2 \lesssim 
    s_0^{(2\alpha-1)/(2\alpha+1)}n^{- 2\alpha/(2\alpha+1)}M_{\alpha,1,n}^{2/(2\alpha+1)} 
    + s_0(\log p)/n. 
\eel
{\rm (v)} When $q = 1$, $q_0=0$, $M_{0,1,n}\vee M_{\alpha,\infty,n}=O(1)$ 
and $\#\{j: \bff^*_j\neq {\bf 0}\} \le s_0$ in \eqref{th-2-4}, 
\bel{cor-1-5}
    \|\hbf - \bff^*\|_{2,n}^2 \lesssim 
    s_0^{1/(2\alpha)} n^{-1/2}
    + s_0(\log p)/n. 
\eel
\end{corollary}

We note that $n^{-\gamma/2}\le \{(\log p)/n\}^{1-q/2}$ 
when $\gamma = 2-q$ and $p\ge 3$ in part (ii) above, and omit the proof of Corollary \ref{cor-1}. 

Compared with the existing results discussed in Section 1, 
\eqref{cor-1-1} 
is of the same form as the oracle inequality 
\eqref{KY} in \cite{koltchinskii2010sparsity} and \cite{raskutti2012minimax}, 
\eqref{cor-1-2} 
is of the same form as the oracle inequality \eqref{YZ} in \cite{yuan2015minimax}, 
\eqref{cor-1-3} 
is of the same form as the oracle inequality \eqref{TZ} in \cite{tan2019doubly}, 
and \eqref{cor-1-4} 
is of the same form as the oracle inequality in 
\cite{suzuki2012fast}. 
As $\|\sum_{j=1}^p f_j^*(x_j)\|_\infty=\sum_{j=1}^p \|f_j^*(x_j)\|_\infty$ when the maximum 
is taken over $(x_1,\ldots,x_p)^\top \in [0,1]^p$, 
\eqref{cor-1-5} 
improves upon the oracle inequality \eqref{RWY} of \cite{raskutti2012minimax} 
by a logarithmic factor in the first term. 
Thus, \eqref{th-2-1} unifies all the above results with a single oracle inequality, 
achieves rate minimaxity in these special cases as discussed in Section 1 and 
elucidated in the referenced papers, 
and provides new error bounds for $0<\alpha_0<\alpha$ and $q\neq 0$. 
More important, while the estimators in the above existing results all use tuning 
parameters depending on the smoothness index $\alpha$ of the underlying component 
functions $f_j$ and at least implicitly on the regularization parameters $q$ and $q_0$, 
the MR-GL procedure does not involve any such tuning parameters. 

Of course, the above comparisons are made under parallel conditions but not always in 
the same setting. 
While Theorem \ref{th-2} concerns deterministic designs 
and involves the prediction factor $C^*_{\rm pred}(\xi,S)$, the above referenced 
papers mainly focus on random designs. 
In Section \ref{3-2}, we will derive $L_2$ oracle inequalities parallel to Theorem \ref{th-2} 
and Corollary \ref{cor-1} under a population compatibility condition with random designs, 
allowing direct comparisons. 

\subsection{Connection to the group Lasso theory}\label{3-1-3} 
To describe the connection of Theorem \ref{th-1} to the theory of group Lasso, we write 
\bel{GL-coef}
\hbbeta = \argmin_{\bb}
\bigg\{\frac{1}{2}\bigg\|\by - \sum_{k=1}^{g^*} \bU_{k}\bb_{k}\bigg\|_{2,n} + 
A_0\sum_{k=1}^{g^*}  \lam_{k} \|\bU_{k}\bb_{k}\|_{2,n}\bigg\} 
\eel
as the group Lasso in the typical setting where $g^*$ is the total number of groups, 
$\bU_k\in\R^{n\times d_{k}}$ 
is composed of design vectors in the $k$-th group 
$\hbox{\rm GR}_k$, $d_{k}=|\hbox{\rm GR}_k|$ and 
$\bb=(\bb_1^\top,\ldots,\bb^\top_{g^*})^\top$ with $\bb_j\in \R^{d_{k}}$. 
This matches \eqref{MR-GL-coef} with the label change $(j,k)\to k$ and $g^*=|\scrK^*|$. 
Let $\bU=(\bU_1,\ldots,\bU_{g^*})\in \R^{n\times d^*}$ with $d^*=\sum_{k=1}^{g^*}d_{k}$. 

Let $\by \sim N(\bff^*,\sigma^2\bI_{n\times n})$ and 
$\lam_{k} = \sigma n^{-1/2}(d_{k}^{1/2}+t_k)$. 
Theorem~\ref{th-1} asserts that 
for any $\bff^*\in \R^n$, $\bbetabar = (\bbetabar_1^\top,\ldots,\bbetabar^\top_{g^*})^\top$ 
and $S\subset \{1,\ldots,g^*\}$, 
\bel{th-1-1a}
&& \big\|\bU\hbbeta - \bU\bbetabar\big\|_{2,n}^2 + \big\|\bU\hbbeta - \bff^*\big\|_{2,n}^2
\\ \nonumber 
&\le& 2\big\|\bU\bbetabar - \bff^*\big\|_{2,n}^2 
+ 8A_0\sum_{k\in S^c}\lam_{k}\|\bU_{k}\bbetabar_{k}\|_{2,n}
+ 4(A_0+1)^2C_{\rm pred}(\xi,S)\sum_{k\in S}\lam_{k}^2
\eel
with at least probability $1-\sum_{k=1}^{g^*} \P\{N(0,1) > t_k\}$, 
where $\xi = (A_0+1)/(A_0-1)$ and 
\bel{C-pred-1a}
C_{\rm pred}(\xi,S) 
= \sup_{\bb: \|\bU\bb\|_{2,n}>0} \frac{\big\{\sum_{k\in S}\lam_{k}\|\bU_{k}\bb_{k}\|_{2,n} 
- \xi^{-1}\sum_{k\in S^c}\lam_{k}\|\bU_{k}\bb_{k}\|_{2,n}\big\}_+^2}
{\|\sum_{k=1}^{g^*} \bU_{k}\bb_{k}\|_{2,n}^2\sum_{k\in S}\lam_{k}^2}. 
\eel
Theorem \ref{th-1} further suggests the choice $S=\{k: \|\bU_{k}\bbetabar_{k}\|_{2,n}\ge A_0\lam_{k}\}$ 
to achieve 
\bel{th-1-3a}
&& \big\|\bU\hbbeta - \bU\bbetabar\big\|_{2,n}^2 + \big\|\bU\hbbeta - \bff^*\big\|_{2,n}^2
\\ \nonumber 
&\le& 2\big\|\bU\bbetabar - \bff^*\big\|_{2,n}^2 
+ C_{\rm pred}^*(\xi,S)\sum_{k=1}^{g^*}\lam_{k}^2\wedge(\lam_{k}\|\bU_{k}\bbetabar_{k}\|_{2,n})
\eel
with $C^*_{\rm pred}(\xi,S)=\max\big\{8A_0^2, 4(A_0+1)^2C_{\rm pred}(\xi,S)\big\}$. 
By allowing arbitrary $\bff^*$, $\bbetabar$ and $S$, 
the oracle inequalities \eqref{th-1-1a} and (\ref{th-1-3a}) 
improve upon the more familiar existing ones 
under the hard group sparsity $S = \{k: \|\bbetabar_k\|_2 > 0\}$ with $\bff^*=\bU\bbetabar$. 
Under the hard group sparsity, 
\eqref{th-1-1a} gives the squared error rate 
$\sum_{k\in S}\lam_k^2 \asymp \sum_{k\in S}(d_{k}+t_k^2)/n$ 
to exploit the benefit of the group sparsity \citep{huang2010benefit}. 

Let $\scrC(\xi,S) = \{\bb: \sum_{S^c}\lam_k\|\bU_{k}\bb_{k}\|_{2,n} \le 
\xi \sum_{S}\lam_k\|\bU_{k}\bb_{k}\|_{2,n}\}$. 
The prediction factor $C_{\rm pred}(\xi,S)$ 
is closely related to the groupwise compatibility coefficient (CC)
\bel{GCC}
    \kappa(\xi,S) = \inf\bigg\{\frac{\|\sum_{k=1}^{g^*}\bU_{k}\bb_{k}\|_{2,n}\|\lam_S\|_2}
    {\sum_{k\in S}\lam_{k}\|\bU_{k}\bb_{k}\|_{2,n}}: \bb\in \scrC(\xi,S) \bigg\}. 
\eel
with $\|\lam_S\|_2=\big(\sum_{k\in S}\lam_{k}^2\big)^{1/2}$. 
This quantity reduces to the CC for the Lasso 
\citep{van2009conditions} when $d_{k}=1$. 
It is clear by definition that $C_{\rm pred}(\xi,S)\le 1/\kappa^2(\xi,S)$. 
\cite{buhlmann2011statistics} defined the groupwise CC as 
\bel{CC}
    \kappa_0(\xi,S) 
    = \inf\bigg\{\frac{\|\sum_{k=1}^{g^*}\bU_{k}\bb_{k}\|_{2,n}\big(\sum_{k\in S}d_{k}\big)^{1/2}}
    {\sum_{k\in S}d_{k}^{1/2}\|\bb_{k}\|_2}: \bb\in \scrC_0(\xi,S) \bigg\} 
\eel
with $\scrC_0(\xi,S) = \{\bb: \sum_{S^c}d_{k}^{1/2}\|\bb_{k}\|_{2} \le 
\xi \sum_{S}d_{k}^{1/2}\|\bb_{k}\|_{2}\}$ to match group Lasso penalties satisfying  
$\lam_k \propto \sqrt{d_k}$. 
The two versions of the groupwise CC are equivalent up to a constant 
in the full nonparametric AM with 
$2^{k_*-1} < 2\log(p/\eps) \le 2^{k_*}$ as in Theorem~\ref{th-2} 
under mild side conditions. Specifically, 
\bel{compare-CC}
C_{\rm pred}(\xi,S)\le 1/\kappa^2(\xi,S) 
\le \{2 c^*/\kappa_0(\xi_0,S_0)\}^2 
\eel
when $c_*\|\bb_{k}\|_2^2\le \|\bU_{k}\bb_{k}\|_{2,n}^2\le c^*\|\bb_{k}\|_2^2$ 
and $\sqrt{2d_{k}} \le \lam_{k}n^{1/2}/\sigma \le 2\sqrt{2 d_{k}}$ for all $\bb_k\in\R^{d_k}$ 
and $1\le k\le g^*$, 
$S=\{{k}: \|\bU_{k}\bbetabar_{k}\|_{2,n} \ge A_0\lam_{k}\}$, 
$S_0 =\{{k}: c^*\|\bbetabar_{k}\|_{2} \ge A_0\lam_{k}\}$ and $\xi_0= 2\xi c^*/c_*$. 
Inequality \eqref{compare-CC} is a consequence of 
$S\subseteq S_0$ and $\scrC(\xi,S)\subseteq\scrC_0(\xi_0,S_0)$ 
and the monotonicity properties of the two versions of the groupwise CC. 
A similar strategy will be used in our study of CC under random design in Section \ref{3-2}.

\section{Random Designs}\label{3-2} 

In the random design setting, out-of-sample squared prediction error $\|\fhat - f^*\|_{L_2}^2$ will be used 
in our analysis to evaluate the performance of the MR-GL estimator in (\ref{MR-GL-pred}), 
in addition to the in-sample squared prediction error $\|\hbf-\bff^*\|_{2,n}^2$ considered in Section \ref{3-1}. 
Moreover, instead of the empirical groupwise compatibility condition, we will impose a theoretical groupwise 
compatibility condition and prove that the empirical groupwise CC in \eqref{GCC} can be bounded from blow by its population version up to a constant factor under the sample size condition 
$n\gg s (\log s)^2(\log {d^*})(\log n)$.


\subsection{Equivalence between the empirical and population conditions on the design}\label{3-2-1} 
The groupwise compatibility condition is closely related to the groupwise Restricted Eigenvalue (RE) condition 
in the sense that the groupwise CC is always no smaller than the groupwise RE. 
They characterize similar desirable properties of the design matrix 
which ensure the performance of regularized least squares and related methods. 
In general, the RE is aimed at bounding the $\ell_2$ error rate for estimating the coefficients 
and the CC is aimed at the prediction error. 

For $\ell_1$ regularized LSE, error bounds for the Lasso and Dantzig selector 
were established in \cite{bickel2009simultaneous}, \cite{koltchinskii2009dantzig} 
and \cite{van2009conditions} among many others. 
\cite{van2009conditions} called CC the $\ell_1$ RE. 
Great effort has been devoted to establishing the empirical RE-type conditions under more interpretable conditions. 
For instance, \cite{bickel2009simultaneous}, \cite{van2009conditions}, \cite{zhang2009some} and \cite{ye2010rate} 
provided lower bounds for the RE and CC in terms of the lower and upper sparse eigenvalues of the 
empirical Gram matrix or its population version. 
\cite{raskutti2010restricted} and \cite{rudelson2012reconstruction} proved that the RE condition is guaranteed by 
its population version without imposing a condition on the upper sparse eigenvalue condition on the design. 
Based on their results, the RE condition is understood to be of a weaker form than 
the restricted isometry property \citep{candes2005decoding, candes2007dantzig} 
and the sparse Riesz condition \citep{zhang2008sparsity}. 
Additionally, \cite{lecue2014sparse} and \cite{van2014higher} proved the empirical RE-type conditions 
respectively under a high-order moment condition and a higher order isotropy condition. 

As we have mentioned in Section \ref{3-1-3}, the prediction factor used in the oracle inequalities in 
Theorems \ref{th-1} and \ref{th-2}
is bounded by the reciprocal of the corresponding 
groupwise CC as defined in \eqref{GCC}. 

For the group Lasso, groupwise RE and compatibility conditions have been used 
to derive oracle inequalities  
\citep{nardi2008asymptotic, huang2010benefit, lounici2011oracle, buhlmann2011statistics, negahban2012unified, mitra2016benefit} 
but little has been done on the groupwise RE-type conditions. 
It is understood that the groupwise RE and CC can be bounded by the lower and upper sparse eigenvalues 
\citep{mitra2016benefit}. However, it is unclear whether the groupwise RE-type conditions are guaranteed by 
its population versions and if so the sample size required. 
Here we carry out a systematic study of the groupwise RE-type conditions 
for uniformly bounded design variables by extending the analysis in \cite{rudelson2008sparse} 
from the $\ell_1$ regularization to general group regularization. 
The core of our analysis is the following theorem which asserts that the sample RE 
in a general form is bounded from below by its population version up to a constant factor. 

Our general result concerns a pair of norms $\|\cdot\|^*_1$ and $\|\cdot\|^*_2$, respectively 
related to the group regularization and the loss function of interest, such that 
\bel{norms}
\|\bb_k\|_1\le \|\bb_k\|^*_1,\quad \sum_{k\in S_0}\|\bb_k\|^*_1 \le s^{1/2}\|\bb\|^*_2, 
\eel
for all $\bb$ satisfying $\sum_{k\in S_0^c} \|\bb_k\|^*_1 < \xi_0 \sum_{k\in S_0} \|\bb_k\|^*_1$ and a constant $s>0$, where $S_0$ is a deterministic nonempty subset of $\{1,\ldots,g^*\}$. 
Recall that we group elements of vectors $\bb\in\R^{d^*}$ by writing 
$\bb=(\bb_1^\top,\ldots,\bb_{g^*}^\top)^\top$ 
with $\bb_k\in \R^{d_{k}}$ and $d^*=\sum_{k=1}^{g^*} d_k$. 
For $\xi_0>0$ and the norms and $S_0$ in \eqref{norms}, 
define the general deterministic cone as 
$$
\scrC_0(\xi_0,S_0) = \scrC_0(\xi_0,S_0;\|\cdot\|^*_1) = \bigg\{\bb \in \R^{d^*}: 
\sum_{k\in S_0^c} \|\bb_k\|^*_1 < \xi_0 \sum_{k\in S_0} \|\bb_k\|^*_1 \bigg\} 
$$
and the corresponding generalized groupwise RE and its population version by 
\bel{RE}
\RE_0(\xi_0,S_0) = \inf_{\bb\in\scrC_0(\xi_0,S_0;\|\cdot\|^*_1)}\frac{\|\bU\bb\|_{2,n}}
{\|\bb\|^*_2},\quad 
\REbar_0(\xi_0,S_0) = \inf_{\bb\in\scrC_0(\xi_0,S_0;\|\cdot\|^*_1)}\frac{\|\bU\bb\|_{L_2,n}}{\|\bb\|^*_2}, 
\eel
where $\bU$ is the design matrix and $\|\bU\bb\|_{L_2,n}=(\E[\|\bU\bb\|_{2,n}^2])^{1/2}$.  
This includes as special cases the groupwise RE in \cite{lounici2011oracle} 
with $\|\bb_k\|^*_1 = (\lam_kn^{1/2}/\sigma)\|\bb_k\|_2$ and 
$\|\bb\|^*_2=\max_{|T|=|S_0|}(\sum_{k\in T}\|\bb_k\|_2^2)^{1/2}$ 
and the groupwise CC in \eqref{CC}.  
We may also take $\|\bb\|^*_2 =\|\bb\|_2$ for groupwise RE.  
The following theorem provides sufficient conditions under which the generalized  
groupwise RE in \eqref{RE} is bounded from below by its population version. 

\begin{theorem}\label{th-4}
Let $\RE_0(\xi_0,S_0)$ and $\REbar_0(\xi_0,S_0)$ be as in \eqref{RE} with a random matrix 
$\bU \in \mathbb{R}^{n\times {d^*}}$ and  
norms $\|\cdot\|^*_1$ and $\|\cdot\|^*_2$ satisfying \eqref{norms}. Suppose $\bU$ has independent rows and 
is uniformly bounded, $\|\bU\|_{\max}\le L_0$ for some constant $L_0$. 
Then, 
\bel{th-4-1}
\E\Big[\big\{1 - \RE_0^2(\xi_0,S_0)/\REbar_0^2(\xi_0,S_0)\big\}_+\Big]
\le \eta,
\eel
with $\eta = C_0L_0(s_1/n)^{1/2}(\log s_1)\sqrt{(\log d^*)(\log n)}$,  
where $s_1$ is a constant satisfying $s_1 \ge (1+\xi_0)^2s/\REbar_0^2(\xi_0,S_0)$ 
with the $s$ in \eqref{norms}, 
$d^*=\sum_{k=1}^{g^*}d_k$,
and $C_0$ is a numerical constant. Moreover, 
for constants $c_0$ and $\eps_1$ in $(0,1)$, 
\bel{th-4-2}
\P\left\{
\begin{matrix}
\big|\|\bU\bb\|_{2,n} - 1\big| \le c_0\ \forall \bb\in\scrC_0(\xi_0,S_0), \|\bU\bb\|_{L_2,n}^2=1
\cr \big|\|\bU\bb\|_{2,n} - 1\big| \le c_0\ \forall \|\bb\|_1\le s_1^{1/2}, \|\bU\bb\|_{L_2,n}^2=1
\cr \RE_0^2(\xi_0,S_0) \ge (1-c_0)\REbar_0^2(\xi_0,S_0)
\end{matrix} \right\} 
 \ge 1 -\eps_1   
\eel
when $5e^{-c_0/\eta} \le \eps_1$ and $4\pi s_1L_0/n\le \eta$. 
\end{theorem} 

The uniform boundedness condition $\|\bU\|_{\max}\le L_0$ is natural in the study of 
the nonparametric AM when the Fourier basis functions are used to build the MR-GL.  
For the CC in \eqref{CC}, $\|\bb_k\|^*_1= d_k^{1/2}\|\bb_k\|_2$,  
$\|\bb_k\|^*_2=\sum_{k\in S_0}d_k^{1/2}\|\bb_k\|_2/s^{1/2}$, $s = \sum_{k\in S_0}d_k$ 
is the total number of variables for the groups with $k\in S_0$, 
and Theorem \ref{th-4} yields 
\bel{RE-bd}
\P\bigg\{\kappa_0(\xi_0,S_0) \ge \sqrt{1-c_0}\,\inf_{\bb\in\scrC_0(\xi_0,S_0)}\frac{\|\bU\bb\|_{L_2,n}s^{1/2}}
{\sum_{k\in S_0}d_k^{1/2}\|\bb_k\|_2}\bigg\} = 1+o(1)
\eel
when $(s/n)(\log s)^2(\log d^*)(\log n) = o(1)$. 
As $s/n\le 1$ is required to have $\kappa_0(\xi_0,S_0)>0$, 
the sample size requirement is optimal up to a logarithmic factor.  
Theorem \ref{th-4} also yields a similar lower bound for the RE 
with the $\ell_2$ loss 
\bel{CC-bd}
\P\bigg\{\inf_{\bb\in\scrC_0(\xi_0,S_0)}\frac{\|\bU\bb\|_{2,n}}{\|\bb\|_2} 
\ge \sqrt{1-c_0}\inf_{\bb\in\scrC_0(\xi_0,S_0)}\frac{\|\bU\bb\|_{L_2,n}}
{\|\bb\|_2} \bigg\} = 1+o(1) 
\eel
under the same sample size condition for the same cone as in \eqref{RE-bd}. 

To bound the CC \eqref{GCC} more directly associated with the prediction factor in \eqref{C-pred-1}, 
we still need to deal with random cones involving the empirical norms $\|\bU_k\bb_k\|_{2,n}$ 
and possibly random $S\subset \{1,\ldots,g^*\}$ instead of the deterministic $S_0$.  
To this end, we provide in the following lemma the equivalence between the empirical norm 
$\|\bU_k\bb_k\|_{2,n}$
and its population version $\|\bU_k\bb_k\|_{L_2,n}$ with an application of 
the non-commutative Bernstein inequality \citep{tropp2012user}.  

\begin{lemma}\label{lemma-1} 
Let $\bU_{k} \in \R^{n\times d_{k}}$ be random matrices with independent rows $\br_k^{i}$ 
satisfying $\P\big\{\|\br_k^{i}\|_2 \leq L_0d_{k}^{1/2}, \forall i\leq n \big\} =1$ for some 
constant $L_0>0$. 
Let $\nu_{+,k} = \max_{\|\bb_k\|_2=1}\E[\|\bU_k\bb_k\|_{2,n}^2]$
and $\nu_{-,k} = \min_{\|\bb_k\|_2=1}\E[\|\bU_k\bb_k\|_{2,n}^2]$.  
Then, 
\begin{align}
    \P\bigg\{\max_{1\leq k\leq g^*}\max_{\|\bb_k\|_2=1} 
\Big| \|\bU_k\bb_k\|_{2,n}^2 - \E[\|\bU_k\bb_k\|_{2,n}^2]\Big| > c_0 \bigg\} \leq \eps_1   
\label{lm-1-1}
\end{align}
when $\sum_{k=1}^{g^*} 2d_k \exp\big[ - nc_0^2\big/
\big\{2d_{k}L_0^2(\nu_{+,k}+ c_0/3)\big\}\big]\le \eps_1$. 
Moreover,  
\begin{align}
    \P\bigg\{\max_{1\leq k\leq g^*}\max_{\E[\|\bU_k\bb_k\|_{2,n}^2]> 0}
\bigg| \frac{\|\bU_k\bb_k\|_{2,n}^2}{\E[\|\bU_k\bb_k\|_{2,n}^2]} - 1\bigg| > c_0 \bigg\} 
\leq \eps_1 
\label{lm-1-2}
\end{align}
when $\sum_{k=1}^{g^*} 2d_k \exp\big[ - n c_0^2/\big\{2d_{k}(L_0^2/\nu_{-,k})
(1+ c_0/3)\big\}\big]\le \eps_1$. 
\end{lemma}

In Lemma \ref{lemma-1}, \eqref{lm-1-1} states that with probability at least $1-\eps_1$,  
the squared empirical norm $\|\bU_{k}\bb_{k}\|_{2,n}^2$ and its population version 
$\|\bU_{k}\bb_{k}\|_{L_2,n}^2 =\E[\|\bU_k\bb_k\|_{2,n}^2]$ differ by at most a small fraction 
of $\|\bb_k\|_2^2$ simultaneously for all $\bb_k$ and $k$, and \eqref{lm-1-2} implies 
\bel{lm-1-3}
    \P \big\{ c_- \|\bb_{k}\|_2^2 \leq \|\bU_{k}\bb_{k}\|_{2,n}^2 
    \leq c_+ \|\bb_{k}\|_2^2\ \forall\ \bb_k, k\le g^* \big\} \geq 1-\eps_1 
\eel
with $c_- = \nu_-(1- c_0)$ and $c_+ \le \nu_+(1+ c_0)$.

With Theorem \ref{th-4} and Lemma \ref{lemma-1} we are ready to study 
the (empirical) groupwise CC. 
Given $\xi_0>0$ and a deterministic subset $S_0$ of $\{1,\ldots,g^*\}$, define 
\bel{pop-CC-GL}
\kappabar(\xi_0,S_0) = \inf_{\bb\in\scrC_0(\xi_0,S_0)}
\frac{\|\bU\bb\|_{L_2,n}\|\lam_{S_0}\|_2}{\sum_{k\in S_0}\lam_k\|\bU_k\bb_k\|_{L_2,n}}
\eel
as the population version of the groupwise CC in \eqref{GCC}, 
with $\|\lam_{S_0}\|_2=(\sum_{k\in S_0}\lam_k^2)^{1/2}$ and 
$\scrC_0(\xi_0,S_0)=\big\{\bb: \sum_{k\in S_0^c}\lam_k\|\bU_k\bb_k\|_{L_2,n}
\le\xi_0\sum_{k\in S_0}\lam_k\|\bU_k\bb_k\|_{L_2,n}\big\}$. 

\begin{corollary}\label{cor-2}
Let $\bU$ be a random matrix with independent rows and 
$\|\bU\|_{\max}\le L_0$ for some constant $L_0$. 
Let $\xi>0$, $S_0$ be a deterministic subset of $\{1,\ldots,g^*\}$ and 
$\kappa(\xi,S_0)$ as in \eqref{GCC} with $S=S_0$. 
Let $c_0\in (0,1)$, $\xi_0=\xi\sqrt{(1+c_0)/(1-c_0)}$ and $\kappabar(\xi_0,S_0)$ 
be as in \eqref{pop-CC-GL}.  
Let $\nu_-=\min_{1\le k\le g^*}\nu_{-,k}$ with the $\nu_{-,k}$ in Lemma \ref{lemma-1}. 
Then, 
\bel{cor-2-1}
\P\big\{ \kappa(\xi,S_0) \ge \sqrt{(1-c_0)/(1+c_0)} \kappabar(\xi_0,S_0) \big\} \ge 1 - 2\eps_1 
\eel
when $c_0$ and $\eps_1$ satisfy the conditions for \eqref{th-4-2} and \eqref{lm-1-2} 
with $s= \|\lam_{S_0}\|_2^2 n/(\sigma^2\nu_-)$. 
\end{corollary}

This can be seen as follows. 
In the event in \eqref{lm-1-2}, $\scrC(\xi,S_0)\subseteq \scrC_0(\xi_0,S_0)$ so that 
$$
\frac{\kappa(\xi,S_0)}{(1+c_0)^{-1/2}} \ge \inf_{\bb\in\scrC_0(\xi_0,S_0)}
\frac{\|\bU\bb\|_{2,n}\|\lam_{S_0}\|_2}{\sum_{k\in S_0}\lam_k\|\bU_k\bb_k\|_{L_2,n}}
= \inf_{\bb\in\scrC_0(\xi_0,S_0)}
\frac{\|\bU\bb\|_{2,n}\|w_{S_0}\|_2}{\sum_{k\in S_0} w_k\|\bU_k\bb_k\|_{L_2,n}}
$$ 
with $w_k= \lam_kn^{1/2}/(\sigma \nu_-^{1/2})$ to guarantee $\|\bb_k\|_1\le \|\bb_k\|^*_1$ 
in the application of \eqref{th-4-2} with $\|\bb_k\|^*_1 = w_k\|\bU_k\bb_k\|_{L_2,n}$ 
and $\|\bb\|^*_2=\sum_{k\in S_0}w_k\|\bU_k\bb_k\|_{L_2,n}/\|w_{S_0}\|_2$. 
Here the scaling of $w_k$ only impacts the bound through $s= \|\lam_{S_0}\|_2^2 n/(\sigma^2\nu_-)$   
in the conditions for \eqref{th-4-2}. 

\subsection{Out-of-sample error bounds in nonparametric AM}\label{3-2-2} 
In this subsection, we derive the $L_2$ prediction error bound for the MR-GL 
in the nonparametric AM. In this model, a response vector $\by = (y_1,\ldots,y_n)^\top$ 
and a random design matrix $\bX=(x_{i,j})_{n\times p}$ are observed such that 
\bel{LM-2}
y_i = f^*(x_{i,1},\ldots,x_{i,p})+\veps_i,\  f^* = f^*(x_1,\ldots,x_p) = \textsum_{j=1}^p 
f^*_j(x_j), 
\eel
where $\veps_i$ are iid $N(0,\sigma^2)$ variables independent of $\bX$. 
We write the components $f^*_j(\cdot)$ in the following multi-resolution expansion 
in a uniformly bounded basis,  
\bel{LM-2a}
f^*_j(x) = \sum_{k=k_*}^\infty f^*_{j,k}(x),\ 
f^*_{j,k}(x) = \sum_{\ell=1}^{2^{k_*\vee(k-1)}} u_{j,k,\ell}(x) \beta^*_{j,k,\ell},\ \|u_{j,k,\ell}\|_\infty\le L_0, 
\eel
for some fixed constant $L_0$. 
For abbreviation, we write $f^*_j = f^*_j(x_j)$ and $f^*_{j,k} = f^*_{j,k}(x_j)$ 
as functions of $(x_1,\ldots,x_p)^\top$ depending on $x_j$ only.  

The MR-GL estimator $\hbbeta$ in \eqref{MR-GL-coef} yields an estimator of the regression function $f^*$, 
\bel{MR-GL-fun}
\fhat = \sum_{j=1}^p\sum_{k=k^*}^{k^*} \fhat_{j,k},\quad \fhat = \fhat(x_1,\ldots,x_p),\quad 
\fhat_{j,k} = \fhat_{j,k}(x_j), 
\eel 
with $\fhat_{j,k}(x_j) = \sum_{\ell=1}^{2^{k_*\vee(k-1)}} u_{j,k,\ell}(x_j) \hbeta_{j,k,\ell}$ 
through the basis functions $u_{j,k,\ell}(\cdot)$ in \eqref{LM-2a}. 

For deterministic or random functions $h(x_1,\ldots,x_p)$, the $L_2$ norm is defined as 
\bel{G_i}
\|h\|_{L_2} = \|\bh\|_{L_2,n} 
= \bigg(\int_{\R^p}h^2(x_1,\ldots,x_p)G(dx_1,\ldots,dx_p)\bigg)^{1/2}
\eel
with $G(x_1,\ldots,x_p) = n^{-1}\sum_{i=1}^n G_i(x_1,\ldots,x_p)$, 
where $G_i$ is the joint distribution of the $i$-th row of the design matrix $\bX$. 
The loss function 
$\|\fhat - f^*\|_{L_2}^2$ 
measures the out-of-sample prediction error 
as it is the expected difference between the predictions by the estimated $\fhat$ 
and true $f^*$ at an out-of-sample random point $\bx\in\R^p$ with joint distribution $G$. 
In \eqref{G_i} $\bh\in\R^n$ is understood as the realization of $h(\cdot)$ in the sample 
and its elements are given by $h(x_{i,1},\ldots,x_{i,p})$ with the $x_{i,j}$ in \eqref{LM-2a}. 

For $\bU\in \R^{n\times d^*}$ and $\bU_{j,k}\in \R^{n\times 2^{k_*\vee(k-1)}}$ in \eqref{U_{j,k}} with
the basis functions $u_{j,k,\ell}(\cdot)$ in \eqref{LM-2a} and $d^*=2^{k^*}p$ 
and the corresponding vectors $\bb\in \R^{d^*}$ and $\bb_{j,k}\in \R^{2^{k_*\vee(k-1)}}$ 
with elements $b_{j,k,\ell}\in \R$, 
the vector version of the norm is useful as it gives 
$\|\bU\bb\|_{L_2,n}$ and $\|\bU_{j,k}\bb_{j,k}\|_{L_2,n}$ the meaning of $\|\bh\|_{L_2,n}$ and $\|\bh_{j,k}\|_{L_2,n}$ with 
$$
h(x_1,\ldots,x_p) = \sum_{j=1}^p \sum_{k=k_*}^{k^*} h_{j,k}(x_j),\quad 
h_{j,k}(x_j) = \sum_{\ell=1}^{2^{k_*\vee(k-1)}} u_{j,k,\ell}(x_j) b_{j,k,\ell}, 
$$ 
for random $\bb$. This saves explicit definition of the function $h$ when the vector $\bb$ is explicitly given. 
For deterministic $\bb$, $\|\bU\bb\|_{L_2,n}^2 = \E\big[\|\bU\bb\|_{2,n}^2\big]$ 
and $\|\bU_{j,k}\bb_{j,k}\|_{L_2,n}^2 = \E\big[\|\bU_{j,k}\bb_{j,k}\|_{2,n}^2\big]$,   
in agreement with the notation in Subsection \ref{3-2-1}. 

Parallel to \eqref{Sobolev-norm} and \eqref{M_alpha_n}, we define for $\alpha >0$ the Sobolev-type norms 
\bel{Sobolev-norm-coef}
\|f_j^*\|_{\alpha,2}^2 = \sum_{k=k_*+1}^\infty 2^{2\alpha k}\|\bbeta^*_{j,k}\|_{2}^2,\quad   
\|f_j^*\|_{\Sobolev,\alpha}^2 = \|\bbeta^*_{j,k_*}\|_{2}^2+\|f_j^*\|_{\alpha,2}^2, 
\eel
with $\bbeta^*_{j,k}=(\beta^*_{j,k,\ell}, \ell\le 2^{k^*\vee(k-1)})^\top$ 
for the coefficients $\beta^*_{j,k,\ell}$ in \eqref{LM-2a}, 
and we measure the overall complexity of the true $f^* = \sum_{j=1}^p f^*_j$ by 
\bel{M_alpha}
M_{\alpha,q}^q=\sum_{j=1}^p\|f^*_j\|_{\alpha,2}^q,\quad 
M_{q}^{q,\BR}=\sum_{j=1}^p \|\beta^*_{j,k_*}\|_{2}^q,
\eel
for all $q \ge 0$, with $q=0$ treated as the limit at $q=0+$. 
For $\alpha > 1/2$, the finiteness of the norm $\|f_j^*\|_{\alpha,2}$ implies 
the uniform absolute convergence of the series in \eqref{LM-2a} as 
\bes
\|f^*_j\|_\infty &\le& \textsum_{k\ge k_*} L_0 2^{(k^*\vee(k-1))/2} \|\bbeta^*_{j,k}\|_\infty
\cr &\le& L_0\big\{2^{k_*/2}+(4(\alpha-1/2)-1)^{-1/2}2^{-(\alpha-1/2)k_*}\big\}\|f_j^*\|_{\Sobolev,\alpha}. 
\ees

To relate the norms in \eqref{Sobolev-norm} and \eqref{Sobolev-norm-coef}, 
define as in Lemma \ref{lemma-1}  
\bel{nu_pm}
\nu_{+} = \max_{(j,k)\in\scrK^*}\max_{\|\bb_k\|_2=1}\|\bU_{j,k}\bb_{j,k}\|_{L_2,n}^2,\quad 
\nu_{-} = \min_{(j,k)\in\scrK^*}\min_{\|\bb_k\|_2=1}\|\bU_{j,k}\bb_{j,k}\|_{L_2,n}^2, 
\eel
with $\scrK^*=\{1,\ldots,p\}\times\{k_*,\ldots,k^*\}$. 
 As $\E[\|\bff^*_j\|_{\alpha,2,n}^2]=\sum_{k>k_*}^\infty 2^{2\alpha k}\|\bU_{j,k}\bbeta^*_{j,k}\|_{L_2,n}^2$, 
\bes
\E\big[\|\bff^*_j\|_{\alpha,2,n}^2\big]\le\nu_+ \|f^*_j\|_{\alpha,2}^2,\ \  
\E\big[\|\bff^*_j\|_{\Sobolev,\alpha,n}^2\big] \le\nu_+ \|f^*_j\|_{\Sobolev,\alpha}^2, 
\ees
with the $\nu_{+}$ in \eqref{nu_pm}, 
so that the norms in \eqref{Sobolev-norm} can be bounded 
by those in \eqref{Sobolev-norm-coef}. 
Similarly, the complexity measures in \eqref{M_alpha_n} can be bounded by their population version 
\eqref{M_alpha} through 
\bel{M-ratio}
\E\big[ M_{\alpha,q,n}^q\big] \le \nu_+^{q/2} M_{\alpha,q}^q,\ 
\E\big[M_{q,n}^{q,\BR}\big] \le \nu_+^{q/2}M_{q}^{q,\BR},\ q\le 2. 
\eel

For constants $\xi_0>0$ and deterministic $S_0\subset \scrK^*$, define the population CC as 
\bel{pop-CC-AM}
\kappabar(\xi_0,S_0) = \inf_{\bb\in\scrC_0(\xi_0,S_0)}
\frac{\|\bU\bb\|_{L_2,n}(\sum_{(j,k)\in S_0}\lam_{j,k}^2)^{1/2}}{\sum_{(j,k)\in S_0}\lam_{j,k}\|\bU_{j,k}\bb_{j,k}\|_{L_2,n}}, 
\eel
where $\scrK^*$ is as in \eqref{nu_pm}, 
$\bU\bb=\sum_{(j,k)\in\scrK^*}\bU_{j,k}\bb_{j,k}$ and 
\bes
\scrC_0(\xi_0,S_0) = \Big\{\bb: \textsum_{(j,k)\in\scrK^*\setminus S_0}\lam_{j,k}\|\bU_{j,k}\bb_{j,k}\|_{L_2,n}
\le \xi_0 \textsum_{(j,k)\in S_0}\lam_{j,k}\|\bU_{j,k}\bb_{j,k}\|_{L_2,n}\Big\}
\ees
To bound the prediction factor \eqref{C-pred-1} by the population CC via Corollary \ref{cor-2} in the error term $B_{S_0}$ in Theorem \ref{th-1}, we need 
\bel{cond-CC}
  5e^{-c_0/\eta} \le \eps_1,\ 4\pi s_1L_0/n\le \eta,\ 
\frac{n c_0^2}{2^{k^*}L_0^2/\nu_{-}} \ge (4/3)\log(np/\eps), 
\eel
for a fixed $c_0\in (0,1)$ and small $\eps_1\in (0,1)$, 
where $L_0$ is as in \eqref{LM-2a}, 
$s_1$ is a constant satisfying $s_1 \ge (1+\xi_0)^2(4s/\nu_-)/\kappabar^2(\xi_0,S_0)$, 
$\eta = C_0L_0(s_1/n)^{1/2}(\log s_1)\sqrt{(\log n)\log(np)}$ 
and $\nu_{-}$ is as in \eqref{nu_pm}, 
with $s = \sum_{(j,k)\in S_0}2^{k^*\vee(k-1)}$ and $C_0$ being the numerical constant in Theorem \ref{th-4}. 
Under \eqref{cond-CC}, $\P\{\Omega_1\}\ge 1-2\eps_1$ with 
\bel{pf-th5-key}
\Omega_1 = \left\{\begin{matrix}
\big|\|\bU\bb\|_{2,n} - 1\big| \le c_0\ \forall \bb\in\scrC(\xi,S_0), \|\bU\bb\|_{L_2,n}^2=1
\cr \big|\|\bU\bb\|_{2,n}/\|\bU\bb\|_{L_2,n} -1\big| \le c_0\ \forall\, \|\bU\bb\|_{L_2,n} \ge 1, \|\bb\|_1\le s_1^{1/2}
\cr \big|\|\bU_{j,k}\bb_{j,k}\|_{2,n}^2 -1\big| \le c_0\ \forall\, \|\bU_{j,k}\bb_{j,k}\|_{L_2,n}=1, (j,k)\in\scrK^*
\cr \kappa(\xi,S_0) \ge \sqrt{(1-c_0)/(1+c_0)} \kappabar(\xi_0,S_0)
\cr M_{\alpha,q,n}^q \le (1+c_0)^q\nu_+^{q/2} M_{\alpha,q}^q,\ 
M_{q,n}^{q,\BR} \le (1+c_0)^q\nu_+^{q/2}M_{q}^{q,\BR}
\end{matrix} \right\} 
\eel
by \eqref{th-4-2} and \eqref{lm-1-2} as in Corollary \ref{cor-2}, 
where $\xi_0=\sqrt{(1+c_0)/(1-c_0)}\xi$, 
\bel{CC-AM}
\kappa(\xi,S_0) = \inf_{\bb\in\scrC(\xi,S_0)}
\frac{\|\bU\bb\|_{2,n}(\sum_{(j,k)\in S_0}\lam_{j,k}^2)^{1/2}}
{\sum_{(j,k)\in S_0}\lam_{j,k}\|\bU_{j,k}\bb_{j,k}\|_{2,n}} 
\eel
and $\scrC(\xi,S_0)=\big\{\bb: \textsum_{(j,k)\in S_0^c}\lam_{j,k}\|\bU_{j,k}\bb_{j,k}\|_{2,n}
\le \xi \textsum_{(j,k)\in S_0}\lam_{j,k}\|\bU_{j,k}\bb_{j,k}\|_{2,n}\big\}$.  
We note that the condition for \eqref{lm-1-2} holds by the last inequality in \eqref{cond-CC} due to 
$d^*=\sum_{(j,k)\in\scrK^*}2^{k_*\vee(k-1)}=2^{k^*}p\le np$. 
Similarly, the condition for \eqref{th-4-2} holds by the first two inequalities in \eqref{cond-CC} 
as $\|\lam_{S_0}\|_2^2 n/(\sigma^2\nu_-) = \sum_{(j,k)\in S_0}(2^{k/2}+\sqrt{2\log(p/\eps)})^2/\nu_-\le 4s/\nu_-$, 
with $2^{k_*-1} \le 2\log(p/\eps)\le 2^{k_*}$ as in Theorem \ref{th-2}. 

We also need to bound the contribution of $f^*_{j,k}$ to the error term $\Deltabar_{S_0}$ in 
Theorem~\ref{th-1} for $(j,k)\in S_0^c$.  
To this end, we will prove 
$\P\big\{\Omega_1\setminus \Omega_2\big\} \le \eps_2$ with
\bel{Omega_2}
\Omega_2 = \Bigg\{\big\|\bff^*-\bffbar\big\|_{2,n}^2
+4A_0\sum_{(j,k)\in\scrK^*\setminus S_0} \lam_{j,k}\|\bU_{j,k}\bbeta^*_{j,k}\big\|_{2,n} \le 
\frac{\sigma^2s_2}{n}\Bigg\}. 
\eel
with the $\scrK^*$ in \eqref{nu_pm}, $\fbar(x_1,\ldots,x_p) = \sum_{(j,k)\in\scrK^*}f^*_{j,k}(x_j)$ 
and 
\bel{s_2}
s_2 &=& \frac{n}{\sigma^2}\bigg[2\|f^*- \fbar\|_{L_2}^2
+ 4A_0(1+c_0)\nu_+^{1/2} \sum_{(j,k)\in \scrK^*\setminus S_0}\lam_{j,k}\|\bbeta^*_{j,k}\|_2
\\ \nonumber && \qquad + 
\min\bigg\{\frac{(C_{\alpha-1/2}L_0M_{\alpha,1})^2|\log \eps_2|}{n^{2(\alpha+\alpha^*)/(2\alpha_*+1)}}, 
\|f^*- \fbar\|_{L_2}^2(1/\eps_2 - 2)_+\bigg\}
\bigg],
\eel
where $\alpha >1/2$, $C_\alpha = 1/(4^\alpha-1)^{1/2}$ 
and $\bff^*$ and $\bffbar$ are the realized vector version of $f^*$ and $\fbar$ respectively. 
To bound the $L_2$ norm above, we notice that 
\bel{approximate3}
\|f^*- \fbar\|_{L_2}^2 = \bigg\| \sum_{j=1}^p\sum_{k=k^*+1}^\infty f^*_{j,k} \bigg\|_{L_2}^2 
\le C_\alpha^2 2^{- 2\alpha k^*}M_{\alpha,1}^2
\le \frac{(C_\alpha M_{\alpha,1})^2}{n^{2\alpha/(2\alpha^*+1)}} 
\eel
by \eqref{approximate2} and \eqref{M-ratio}, 
and we may impose parallel to \eqref{approximate1} the condition that 
\bel{approximate4}
\|f^*- \fbar\|_{L_2}^2
\le C^*_2 \sum_{j=1}^p\sum_{k=k^*+1}^\infty \big\|  f^*_{j,k} \big\|_{L_2}^2
\le C^*_22^{- 2\alpha k^*}M_{\alpha,2}^2
\le \frac{C^*_2 M_{\alpha,2}^2}{n^{2\alpha/(2\alpha^*+1)}} 
\eel
with some constant $C^*_2$. 
As \eqref{approximate4} is imposed only on the ultra high resolution components of the true $f^*$, 
it is in a much weaker form than the parallel side condition $\big\| \sum_{j=1}^p f_{j} \big\|_{L_2}^2 
\lesssim \sum_{j=1}^p \big\|  f_{j} \big\|_{L_2}^2\ \forall f_j=f_j(x_j)$ typically imposed in the literature 
as in the references discussed below Corollary \ref{cor-1}. 

We are now ready to present our main oracle inequality for the nonparametric AM with random design. 

\begin{theorem}\label{th-5}
Let $f^*$ and $f^*_{j,k}$ be as in \eqref{LM-2} and \eqref{LM-2a}, $\fbar$ as in \eqref{Omega_2}, 
and $\fhat$ as in \eqref{MR-GL-fun} 
with the penalty levels $\lam_{j,k} = \sigma(2^{k/2}+\sqrt{2\log(p/\eps)})/n^{1/2}$ 
and $\{A_0, k_*,k^*,\alpha_*,\alpha\}$ as in Theorem \ref{th-2}. 
Let $\bff^*\in \R^n$, $\bffbar\in \R^n$ and $\hbf\in \R^n$ 
be respectively the realizations of $f^*$, $\fbar$ and $\fhat$ in the sample as in Theorem \ref{th-2}. 
Let $\Omega_0$, $\Omega_1$ and $\Omega_2$ be as in 
\eqref{lam-gen}, \eqref{pf-th5-key} and \eqref{Omega_2} respectively. 
Let $\{\xi, \xi_0,c_0,s_1\}$ be as in \eqref{cond-CC} and \eqref{pf-th5-key}, 
$\scrK^*$ as in \eqref{nu_pm}, $s_2$ as in \eqref{s_2}   
and $\kappabar(\xi_0,S_0)$ the population CC in \eqref{pop-CC-AM} 
with a deterministic $S_0\subset \scrK^*$. \\
(i) Let $\xi=(A_0+1)/(A_0-1)$.
In the event $\Omega_0\cap\Omega_1\cap\Omega_2$, 
\bel{th-5-1}
\big\|\hbf - \bffbar\big\|_{2,n}^2 + \big\|\hbf - \bff^*\big\|_{2,n}^2
\le \frac{4(A_0+1)^2\|\lam_S\|_2^2}{\{(1-c_0)/(1+c_0)\}\kappabar^2(\xi_0,S_0)} +\frac{2\sigma^2s_2}{n}. 
\eel
(ii) Let $s_1\ge s_2$ and $\xi=(3A_0+1)/(A_0-1)$. 
In the event $\Omega_0\cap\Omega_1\cap\Omega_2$, 
\bel{th-5-2}
\big\|\fhat - \fbar\big\|_{L_2}^2+\big\|\fhat - f^*\big\|_{L_2}^2
\le C_{A_0,\nu_-,c_0}\big(\|\lam_S\|_2^2/\kappabar^2(\xi_0,S_0)+\sigma^2s_2/n\big), 
\eel
where $C_{A_0,\nu_-,c_0}$ is a constant depending on $(A_0,\nu_-,c_0)$ only 
with the $\nu_-$ in \eqref{nu_pm}. \\
(iii) Suppose \eqref{LM-2}, \eqref{LM-2a} and \eqref{cond-CC} hold with $\eps_1\in (0,1)$ 
and $s_2\le s_1$. 
Then, 
\bes
\P\{\Omega_0^c\}\le \eps/(2\log(p/\eps))^{1/2},\quad 
\P\{\Omega_1^c\}\le 2\eps_1,\quad \P\{\Omega_1\setminus \Omega_2\} \le \eps_2. 
\ees
(iv) Suppose \eqref{LM-2} and \eqref{LM-2a} hold with a fixed $L_0$, $1/\kappabar(\xi_0,S_0)=O(1)$ and $1/\nu_-=O(1)$ in \eqref{nu_pm}. 
Let $s=\sum_{(j,k)\in S_0}2^{k_*\vee(k-1)}$. 
Suppose $n\gg (s+s_2)(\log(s+s_2))^2(\log n)\log(np)$
and $n \gg 2^{k^*}\log(np)$. Then, 
\bel{th-5-3}
\big\|\fhat - f^*\big\|_{L_2}^2 + \big\|\hbf - \bff^*\big\|_{2,n}^2
\lesssim \sigma^2(s+s_2)/n. 
\eel
(v) Suppose the conditions in (iv) hold with $S_0=\{(j,k)\in\scrK^*: \|\bbeta^*_{j,k}\|_2\ge\lam_{j,k}\}$, 
\eqref{approximate4} holds, $\nu_+=O(1)$ in \eqref{nu_pm} and $\alpha_0\neq 1/q-1/2$. Then,
the oracle inequalities \eqref{th-2-2}, \eqref{th-2-3} and \eqref{th-2-4} all hold 
with $\|\hbf - \bff^*\|_{2,n}^2$ replaced by 
$\big\|\fhat - f^*\big\|_{L_2}^2 + \big\|\hbf - \bff^*\big\|_{2,n}^2$ on the left-hand side and 
$\{M_{\alpha,q,n}, M_{q_0,n}^{q_0,\BR}\}$ replaced by the population 
version $\{M_{\alpha,q},M_{q_0}^{q_0,\BR}\}$ in \eqref{M_alpha} on the right-hand side, 
and with the respective $\{\alpha_0, q, q_0, q_1, q_2, \gamma, \rho\}$ in Theorem \ref{th-2}. 
\end{theorem}

We note that Theorem \ref{th-5} (i) and (ii) are analytical and probability is involved only in 
Theorem \ref{th-5} (iii), (iv) and (v). 
In view of Theorem \ref{th-2} and Corollary \ref{cor-1}, Theorem \ref{th-5} (v) implies the 
following corollary. 
Theorem \ref{th-5} uses quantities $\nu_\pm$ in \eqref{nu_pm} to 
bound the norms $\|\bb_{j,k}\|_2$ and $\|\bU_{j,k}\bb_{j,k}\|_{L_2,n}$ by each other with these constant factors, 
largely due to the use of the $\ell_2$ norm to define the complexity 
measures in \eqref{Sobolev-norm-coef} and \eqref{M_alpha}.

\begin{corollary}\label{cor-3} 
Suppose the conditions of Theorem \ref{th-5} (v) hold with
$\alpha\ge\alpha_*$, $\alpha_0\neq 1/q-1/2$ and $\log(1/\eps) = O(\log p)$. \\
{\rm (i)} If $\#\{j: f^*_j\neq 0\} \le s_0$ and $M_{\alpha,\infty}=O(1)$, then   
\bel{cor-3-1}
    \big\|\fhat - f^*\big\|_{L_2}^2 + \big\|\hbf - \bff^*\big\|_{2,n}^2 \lesssim  
    s_0n^{-2\alpha/(2\alpha+1)} + s_0(\log p)/n. 
\eel
{\rm (ii)} If $M_{\alpha,q}^q\vee M_{q}^{q,\BR}=O(1)$, then
\bel{cor-3-2}
    \big\|\fhat - f^*\big\|_{L_2}^2 + \big\|\hbf - \bff^*\big\|_{2,n}^2 \lesssim 
    n^{- 2\alpha/(2\alpha+1)}
    + \big((\log p)/n\big)^{1-q/2}. 
\eel
{\rm (iii)} If $M_{\alpha,1}\vee M_{0,q}^q\vee M_{q}^{q,\BR}=O(1)$, then
\bel{cor-3-3}
    \big\|\fhat - f^*\big\|_{L_2}^2 + \big\|\hbf - \bff^*\big\|_{2,n}^2 \lesssim 
    n^{- (2-q)\alpha/(2\alpha+1-q)}
    + \big((\log p)/n\big)^{1-q/2}. 
\eel
{\rm (iv)} If $\#\{j: f^*_j\neq 0\} \le s_0$, then 
\bel{cor-3-4}
    \big\|\fhat - f^*\big\|_{L_2}^2 + \big\|\hbf - \bff^*\big\|_{2,n}^2 \lesssim 
    s_0^{(2\alpha-1)/(2\alpha+1)}n^{- 2\alpha/(2\alpha+1)}M_{\alpha,1}^{2/(2\alpha+1)} 
    + s_0(\log p)/n. 
\eel
{\rm (v)} If $M_{0,1}\vee M_{\alpha,\infty}=O(1)$ 
and $\#\{j: f^*_j\neq 0\} \le s_0$, then 
\bel{cor-3-5}
    \big\|\fhat - f^*\big\|_{L_2}^2 + \big\|\hbf - \bff^*\big\|_{2,n}^2 \lesssim 
    s_0^{1/(2\alpha)} n^{-1/2}
    + s_0(\log p)/n. 
\eel
\end{corollary}

As we discussed below Corollary \ref{cor-1}, \eqref{cor-3-1}, \eqref{cor-3-2}, \eqref{cor-3-3}, 
\eqref{cor-3-4} and \eqref{cor-3-5} are directly comparable respectively with 
\eqref{KY} of \cite{koltchinskii2010sparsity} and \cite{raskutti2012minimax}, 
\eqref{YZ} of \cite{yuan2015minimax}, 
\eqref{TZ} of \cite{tan2019doubly}, 
the oracle inequality of \cite{suzuki2012fast}
and \eqref{RWY} of \cite{raskutti2012minimax} 
under population compatibility conditions for random design.

\section{Appendix} 

{\sc Proof of Theorem \ref{th-1}.} 
(i) Applying Lemma 4 of \cite{tan2019doubly} with their $\|\cdot\|_{F,j}=0$ to functional effects $f_{j,k}$ , 
we find that for any $\bbetabar$ and $S\subset \scrK^*$
\bes
&& \frac{1}{2}\big\|\hbf - \bffbar\big\|_{2,n}^2 + \frac{1}{2}\big\|\hbf - \bff^*\big\|_{2,n}^2 
+ (A_0-1)\pen_{S^c}(\hbbeta -\bbetabar) 
\cr &\le& \frac{1}{2}\big\|\bffbar- \bff^*\big\|_{2,n}^2 + 2A_0\pen_{S^c}(\bbetabar) 
+ (A_0+1)\pen_{S}(\hbbeta -\bbetabar) . 
\ees
If $(A_0+1)\pen_{S}(\hbbeta -\bbetabar)  \le (A_0-1)\pen_{S^c}(\hbbeta -\bbetabar)$, 
we have 
\bes
\big\|\hbf - \bffbar\big\|_{2,n}^2 + \big\|\hbf - \bff^*\big\|_{2,n}^2 
\le \big\|\bffbar- \bff^*\big\|_{2,n}^2 + 4A_0\pen_{S^c}(\bbetabar) 
= \Deltabar_S.   
\ees
Otherwise, by the definition of $C_{\rm pred}(\xi,S)$, 
\bes
\big\{\pen_n(\{\bh_{j,k}\};S) - \xi^{-1}\pen_n(\{\bh_{j,k}\};S^c)\big\}_+^2
\le C_{\rm pred}(\xi,S)\|\lam_S\|_2^2 = B_S/(A_0+1)^2
\ees
with $\bh_{j,k}=(\hbf_{j,k}-\bffbar_{j,k})/\|\hbf-\bffbar\|_{2,n}$, so that 
\bes
\frac{1}{2}\big\|\hbf - \bffbar\big\|_{2,n}^2 + \frac{1}{2}\big\|\hbf - \bff^*\big\|_{2,n}^2 
\le \frac{1}{2}\Deltabar_S
+ B_S^{1/2}\big\|\hbf - \bffbar\big\|_{2,n}. 
\ees
This gives $\big\|\hbf - \bff^*\big\|_{2,n}^2\le \Deltabar_S+B_S$. 
We also have $z^2/2\le \Deltabar_S/2+B_S^{1/2}z$ with 
$z = \big(\big\|\hbf - \bffbar\big\|_{2,n}^2 + \big\|\hbf - \bff^*\big\|_{2,n}^2\big)^{1/2}$, 
so that $z^2\le (B_S^{1/2}+\sqrt{B_S+\Deltabar_S})^2\le 4B_S+2\Deltabar_S$. 
For \eqref{th-1-3}, we note that $\|\bffbar_{j,k}\|_{2,n}\le A_0\lam_{j,k}$ on $\scrK^*\setminus S$. 

(ii) When $\by - \bff^*$ has iid $N(0,\sigma^2)$ entries, 
$\|\bP_{j,k}(\by - \bff^*)\|_{2,n}^2n/\sigma^2$ has the chi-square 
distribution with $\rank(\bP_{j,k}) \le d_{j,k}$ degrees of freedom. 
Thus, by \eqref{lam} and the Gaussian 
concentration inequality \citep{borell1975brunn, kwapien1994remark} 
\bes
&& \P\left\{\sup_{(j,k)\in\scrK^*} \|\bP_{j,k}(\by - \bff^*)\|_{2,n}/\lam_{j,k} > 1 \right\} 
\cr &\le& 
\sum_{(j,k)\in\scrK^*}\P\Big\{ N(0,1)  > \lam_{j,k}\sqrt{n}/\sigma - \sqrt{d_{j,k}}\Big\} 
\cr &\le&  \frac{\eps}{\sqrt{4\pi\log(p/\eps)}} + |J_1| \sum_{k=k_*+1}^{k^*} 
\frac{\exp\big[ - \big\{2^{k/2}-2^{(k-1)/2}+\sqrt{2\log(p/\eps)}\big\}^2/2\big]}{\sqrt{4\pi\log(p/\eps)}}
\cr &\le& \frac{\eps}{\sqrt{2\log(p/\eps)}}\bigg\{\frac{1}{\sqrt{2\pi}} 
+ \sum_{k=1}^{\infty}\frac{e^{-(1-2^{-1/2})^22^{k-1}-(1-2^{-1/2})2^{k}}}{\sqrt{2\pi}}\bigg\} 
\cr &\le& \frac{\eps}{\sqrt{2\log(p/\eps)}}
\ees
when $2\log(p/\eps)\ge 1$. In the above, 
$\lam_{j,k}\sqrt{n}/\sigma - \sqrt{d_{j,k}} = \sqrt{2\log(p/\eps)}$ for $j\in J_0$ ($k=k_*$ necessarily) 
and $\lam_{j,k}\sqrt{n}/\sigma - \sqrt{d_{j,k}} =2^{k/2}-2^{(k-1)/2}+\sqrt{2\log(p/\eps)}$ for $j\in J_1$. 
$\hfill\square$

\medskip
{\sc Proof of Proposition \ref{prop-1}.} 
Let $t_0>0$ and 
\bes
\kbar = \begin{cases} k_*, & \gamma = 1, 
\cr \big\{k_*\vee\big\lfloor \log_2\big(\sigma_n^{-(\gamma-1)/(\alpha - 1/2)}/t_0\big)\big\rfloor\big\}
\wedge k^*, & \gamma > 1,
\end{cases}
\ees
with integers $0\le k_*\le k^*$. 
We shall first prove 
\bel{prop-11} 
&& \sum_{k=k_*+1}^{k^*}\lam_{k}\min\Big(\|\bffbar_{j,k}\|_{2,n},\lam_{k}\Big)
\\ \nonumber &\le& \sigma_n^{\gamma}
\big\{t_0^{-(1-q/2-\alpha_0q)_+}J_{1-q/2-\alpha_0q}^{(q)}(k_*,\kbar)\|\bffbar_{j}\|_{\alpha_0,2,n}^q
+ t_0^{\alpha-1/2}J_{\alpha-1/2}^{(1)}(\kbar,k^*)\|\bffbar_j\|_{\alpha,2,n} \big\}. 
\eel

As $\lam_k \le \sigma_n 2^{k/2}$, we have 
\bel{pf-prop-1} 
\sum_{k=k_*+1}^{k^*}\lam_{k}\min\Big(\|\bffbar_{j,k}\|_{2,n},\lam_{k}\Big) 
\le \sum_{k=k_*+1}^{k^*}\sigma_n2^{k/2}\min\Big(\|\bffbar_{j,k}\|_{2,n},\sigma_n2^{k/2}\Big). 
\eel
Let $\tbar = \sigma_n^{-(\gamma-1)/(\alpha - 1/2)}/t_0$ for $\alpha >1/2$ 
and $\tbar=2^{k_*}$ for $\alpha=1/2$.  
We write $\kbar$ as the positive integer satisfying $k_*\le \kbar\le k^*$, 
$2^{\kbar} \le \tbar \vee 2^{k_*}$ and $\tbar \wedge 2^{k^*+1} \le 2^{\kbar +1}$. 
For $1-q/2 -\alpha_0 q \le 0$, we have $\gamma = 2-q$ and 
\bes
&& \sum_{k=k_*+1}^{k^*}\sigma_n2^{k/2}\min\Big(\|\bffbar_{j,k}\|_{2,n},\sigma_n2^{k/2}\Big)
\cr &\le & \sum_{k_* < k \le \kbar }\Big(\sigma_n2^{k/2}\Big)^{2-q} 
2^{-\alpha_0 k q}\|2^{\alpha_0 k}\bffbar_{j,k}\|_{2,n}^q 
+ \sum_{\kbar < k \le k^*}\sigma_n2^{- k(\alpha - 1/2)}\|2^{\alpha k}\bffbar_{j,k}\|_{2,n}
\cr &\le & \sigma_n^{2-q}
\sum_{k_* < k \le \kbar }2^{k(1-q/2 -\alpha_0q)} \|2^{\alpha_0 k}\bffbar_{j,k}\|_{2,n}^q 
+ \sigma_n\sum_{\kbar < k \le k^*}2^{- k(\alpha - 1/2)}\|2^{\alpha k}\bffbar_{j,k}\|_{2,n}
\cr &\le & J_{1-q/2-\alpha_0q}^{(q)}(k_*,\kbar)\sigma_n^{\gamma}\|\bffbar_{j}\|_{\alpha_0,2,n}^q 
+ J_{\alpha-1/2}^{(1)}(\kbar,k^*)\sigma_n2^{- (\kbar +1)(\alpha - 1/2)}\|\bffbar_{j}\|_{\alpha,2,n}. 
\ees
Thus, as $2^{- (\kbar +1)(\alpha - 1/2)}\le \tbar^{1/2-\alpha} 
= \sigma_n^{\gamma - 1}t_0^{\alpha-1/2}$ for $\kbar<k^*$ and $J_{\alpha-1/2}^{(1)}(k^*,k^*)=0$, 
\bel{pf-prop-2}
&& \sum_{k=k_*+1}^{k^*}\sigma_n2^{k/2}\min\Big(\|\bffbar_{j,k}\|_{2,n},\sigma_n2^{k/2}\Big)\\ \nonumber
&\le & J_{1-q/2-\alpha_0q}^{(q)}(k_*,\kbar)
t_0^{-(1-q/2 -\alpha_0q)_+}\sigma_n^{\gamma}\|\bffbar_{j}\|_{\alpha_0,2,n}^q + J_{\alpha-1/2}^{(1)}(\kbar,k^*)t_0^{\alpha-1/2}\sigma_n^{\gamma}\|\bffbar_{j}\|_{\alpha,2,n}. 
\eel
Similarly, for $1-q/2 -\alpha_0q > 0$, we have 
\bes
&& \sum_{k=k_*+1}^{k^*}\sigma_n2^{k/2}\min\Big(\|\bffbar_{j,k}\|_{2,n},\sigma_n2^{k/2}\Big)
\cr &\le & \sigma_n^{2-q}
\sum_{k_* < k \le \kbar }2^{k(1-q/2 -\alpha_0q)} \|2^{\alpha_0 k}\bffbar_{j,k}\|_{2,n}^q 
+ \sum_{\kbar < k \le k^*}\sigma_n2^{- k(\alpha - 1/2)}\|2^{\alpha k}\bffbar_{j,k}\|_{2,n}
\cr &\le & J_{1-q/2-\alpha_0q}^{(q)}(k_*,\kbar)
\sigma_n^{2-q}\tbar^{(1-q/2 -\alpha_0q)_+}\|\bffbar_{j}\|_{\alpha_0,2,n}^q 
+ J_{\alpha-1/2}^{(1)}(\kbar,k^*)\sigma_n\tbar^{1/2-\alpha}\|\bffbar_{j}\|_{\alpha,2,n}. 
\ees
For $\alpha > 1/2$, $(2-q) - (1-q/2-\alpha_0q)_+(\gamma-1)/(\alpha-1/2)=\gamma$ by algebra. 
For $\alpha = 1/2$, $q=1$, $\kbar = k_*$ and $J_{1-q/2-\alpha_0q}^{(q)}(k_*,\kbar)=0$. Thus, by the definition of $\tbar$, $\sigma_n^{2-q}\tbar^{1-q/2 -\alpha_0q} = \sigma_n^\gamma t_0^{-(1-q/2 -\alpha_0q)_+}$, 
so that (\ref{pf-prop-2}) also holds for $1-q/2 -\alpha_0q > 0$. Consequently, (\ref{prop-11}) follows from (\ref{pf-prop-1}) and (\ref{pf-prop-2}). Finally, we obtain (\ref{prop-12}) by minimizing 
\bes
t_0^{-(1-q/2-\alpha_0q)_+}J_{1-q/2-\alpha_0q}^{(q)}(k_*,k^*)\|\bffbar_{j}\|_{\alpha_0,2,n}^q+ t_0^{\alpha-1/2}J_{\alpha-1/2}^{(1)}(k_*,k^*)\|\bffbar_{j}\|_{\alpha,2,n}
\ees
over $t_0>0$. 
Note that (\ref{prop-12}) follows from (\ref{prop-11}) directly when $\alpha = \alpha_0 = 1/2=q/2$. 
$\hfill\square$

\medskip
{\sc Proof of Theorem \ref{th-2}.} 
It follows from \eqref{approximate1} and \eqref{approximate2} that 
$2\|\bff^* - \bffbar\|_{2,n}$, the first term in \eqref{th-1-3}, 
is bounded by the first term on the right-hand side of \eqref{th-2-1}. 
For the second term in \eqref{th-1-3}, Proposition \ref{prop-1} implies that the summation is bounded by 
\bes
&& \sum_{(j,k)\in\scrK^*}\lam_{j,k}^2\wedge\big(\lam_{j,k}\|\bffbar_{j,k}\|_{2,n}\big)
\cr &\le & 4 \sum_{j\in J_1}\sum_{k=k_*+1}^{k^*} (\sigma_n^22^k)
\wedge\big(\sigma_n 2^{k/2}\|\bffbar_{j,k}\|_{2,n}\big)
+ \sum_{j=1}^p \lam_{j,k_*}\big(\lam_{j,k_*} \wedge\|\bffbar_{j,k_*}\|_{2,n}\big)
\cr&\le & 4 \sum_{j\in J_1}\sigma_n^{\gamma}J_{q,\alpha,\alpha_0}(k_*,k^*)
\big\{\|\bff^*_{j}\|_{\alpha_0,2,n}^q\big\}^{1-\rho}\|\bff^*_{j}\|_{\alpha,2,n}^{\rho}
+ \lam_0^{2-q_0} \sum_{j=1}^p w_{j,k_*}^{2-q}\|\bffbar_{j,k_*}\|_{2,n}^q
\cr&\le & 4\sigma_n^\gamma 
J_{q,\alpha,\alpha_0}(k_*,k^*) M_{\alpha_0,q_2,n}^{q(1-\rho)}M_{\alpha,q_1,n}^{\rho}
+ \lam_0^{2-q_0}M_{q_0,n}^{q_0,\BR}
\ees
as $\bffbar_{j,k}=\bff^*_{j,k}$. 
This gives \eqref{th-2-1}. 
When $C^*_{\rm pred}(\xi,S)=O(1)$, \eqref{th-2-2}, \eqref{th-2-3} and \eqref{th-2-4} follows 
from \eqref{th-2-1} by plugging-in the respective values of $(\gamma, q, \rho, q_1, q_2)$. 
For \eqref{th-2-2}, $M_{\alpha,2,n}^2\le s_0M_{\alpha,\infty,n}^2$ 
and $M_{\alpha_0,q_2,n}^{q(1-\rho)}M_{\alpha,q_1,n}^{\rho}
\le (s_0M_{\alpha,\infty,n}^q)^{1-\rho}(s_0M_{\alpha,\infty,n})^{\rho}=s_0M_{\alpha,\infty,n}^{2-\gamma}$. 
For \eqref{th-2-3}, $M_{\alpha,1,n} \le M_{\alpha,q,n}$ and 
$q(1-\rho)+\rho=2-\gamma$ by Remark \ref{remark-1}. 
$\hfill\square$


\medskip
{\sc Proof of Theorem \ref{th-4}}. 
By the scale invariance of the ratios in \eqref{RE}, 
\bes
\RE_0(\xi,S_0) = \inf_{\bb\in\scrC_1(\xi,S_0)}\frac{\|\bU\bb\|_{2,n}}{\|\bb\|^*_2} 
\ge \REbar_0(\xi,S_0)\inf_{\bb\in\scrC_1(\xi,S_0)}\|\bU\bb\|_{2,n} 
\ees
with  $\scrC_1(\xi,S_0) = \big\{\bb: \bb \in \scrC_0(\xi,S_0), \E\big[\|\bU\bb\|_{2,n}^2\big]=1\big\}$. 
For $\bb\in \scrC_1(\xi,S_0)$, 
\bes
\|\bb\|_1\le \sum_{k=1}^{g^*}\|\bb_k\|^*_1 
\le (1+\xi)\sum_{k\in S_0}\|\bb_k\|^*_1
\le (1+\xi)s^{1/2}\|\bb\|^*_2
\le \frac{(1+\xi)s^{1/2}}{\REbar_0(\xi,S_0)}
\le s_1^{1/2}. 
\ees 
Recall that $s_1$ is a constant satisfying $s_1 \ge (1+\xi)^2s/\REbar_0^2(\xi,S_0)$.  
Define 
\bel{M_U}
M_{\bU} = \sup_{\bb\in\scrC_1^*(s_1)}\big| \|\bU\bb\|_{2,n}^2 - 1\big| 
\eel
with $\scrC_1^*(s_1) = \{ \bb\in\R^{d^*}: \|\bb\|_1\le s_1^{1/2}, \E\big[\|\bU\bb\|_{2,n}^2\big]=1\}$.   
We have 
\bel{M_U-2}
\big\{1 - \RE_0^2(\xi,S_0)/\REbar_0^2(\xi,S_0)\big\}_+ \le M_{\bU}. 
\eel

\def\sym{\hbox{\rm\footnotesize sym}}
Let $\br^i$ be the $i$-th row of $\bU$. 
To bound $M_{\bU}$, we write it as 
\bel{M_U-1}
M_{\bU} 
= \sup_{\bb\in\scrC_1^*(s_1)}\bigg|\frac{1}{n}\sum_{i=1}^n 
\Big(\langle \br^i,\bb\rangle^2 - \E\langle \br^i,\bb\rangle^2\Big) \bigg|. 
\eel
As in \cite{rudelson2008sparse}, we shall apply 
the Gaussian symmetrization and Dudley's inequality to bound the expectation $\E[M_{\bU}]$.
By the Gaussian symmetrization 
[Lemma 6.3 and (4.8) in \cite{ledoux1991probability}], 
\bel{Dudley-0} 
 \E\big[M_{\bU}\big] 
\leq \E\big[M_{\bU}^{\sym}\big]
\ \hbox{ with }\ 
 M_{\bU}^{\sym} = \sup\bigg\{\bigg|\frac{\sqrt{2\pi}}{n}\sum_{i=1}^n 
 z_i \langle \br^i,\bb\rangle^2 \bigg|:  \bb\in \scrC_1^*(s_1) \bigg\}, 
\eel
where $\{z_i, i\le n\}$ is a set of independent $N(0,1)$ variables independent of $\bU$. 
Let $\E_{\bU}$ be the conditional expectation given $\bU$ 
and $\mu_{\bU}$ the conditional expectation of $M_{\bU}^{\sym}$. 
By Dudley's inequality [Theorem 11.17 in \cite{ledoux1991probability}],  
\bel{Dudley} 
\mu_{\bU} = \E_{\bU}\big[M_{\bU}^{\sym}\big] 
\le \frac{24\sqrt{2\pi}}{n} 
\int_0^\infty \sqrt{\log N\big( \scrC_1^*(s_1), d(\cdot , \cdot), w\big)} dw 
\eel
where
$N\big(\scrC_1^*(s_1), d(\cdot , \cdot), w\big)$ is the minimal covering number of 
$\scrC_1^*(s_1)$ by balls of radius $w$ in the metric $d(\cdot , \cdot)$ given by 
$$
d(\ba, \bb) = \big\{\textsum_{i=1}^n \big(\langle\br^i,\ba\rangle^2 - \langle\br^i,\bb\rangle^2 \big)^2\big\}^{1/2}.
$$ 
To bound the entropy integral in \eqref{Dudley}, we 
transform the distance into a simpler one as follows. 
As $\max_{\bb\in\scrC_1^*(s_1)}\sum_{i=1}^n \langle\br^i,\bb\rangle^2/n \le M_{\bU}+1$ by \eqref{M_U}, 
for $\{\ba,\bb\}\subset \scrC_1^*(s_1)$
\begin{align*}
    d(\ba, \bb)     
    \leq L_0 \bigg[ \sum_{i=1}^n \langle\br^i,\ba+\bb\rangle^2 \bigg]^{1/2}\|\ba-\bb\|_{\bU} 
    \leq 2L_0 n^{1/2}(M_{\bU}+1)^{1/2}\|\ba-\bb\|_{\bU}, 
\end{align*}
where $\|\bb\|_{\bU} = \max_{i\leq n} |\langle \bb, \br^i \rangle|/L_0$.   
Let $\calB_1^{{d^*}}=\{\bb: \|\bb\|_1\le 1\}$ and 
\bes
J_{\bU} = \int_0^1 \sqrt{\log N\big(\calB_1^{{d^*}}, \|\cdot\|_{\bU}, w\big)} dw
\ees
As $\|\bb\|_1\le s_1^{1/2}$ in $\scrC_1^*(s_1)$ and 
$d(\ba, \bb) \leq 2L_0 n^{1/2}(M_{\bU}+1)^{1/2}\|\ba-\bb\|_{\bU}$, 
\eqref{Dudley} is further bounded by 
\bes
 \mu_{\bU}
 &\leq& \frac{24\sqrt{2\pi}}{n^{1/2}} \big(2L_0 (M_{\bU}+1)^{1/2}\big) 
 \int_0^\infty \sqrt{\log N\big(s_1^{1/2}\calB_1^{{d^*}}, \|\cdot\|_{\bU}, w\big)} dw
 \cr & = & C_0' L_0 (s_1/n)^{1/2}(M_{\bU}+1)^{1/2}J_{\bU}
\ees
with $C_0'=48\sqrt{2\pi}$. 
As $\|J_{\bU}\|_\infty \le C_0''\log(s_1)\sqrt{(\log d^*)\log n}$ by (3.9) of \cite{rudelson2008sparse}, we have 
\bel{new-pf-1}
 \mu_{\bU} \le  (M_{\bU}^{1/2}+1)\eta/4 \le  M_{\bU}/4  + (\eta^2+\eta)/4
\eel
with $\eta = 4C_0'L_0(s_1/n)^{1/2}C_0''\log(s_1)\sqrt{(\log d^*)\log n}$. 
Thus, by Cauchy-Schwarz,
\bes
 \E\big[M_{\bU}\big] \le \E\big[\mu_{\bU}\big] 
 \le (\eta^2+\eta)/3 \le \eta\ \hbox{ when }\ \eta < 1. 
 \ees
We note that \eqref{th-4-1} is trivial when $\eta>1$. 

The tail probability bound can be derived in a similar symmetrization argument. 
Let 
$\sigma_{\bU}^2 = \sup_{\bb\in\scrC_1^*(s_1)}\sum_{i=1}^n 
\langle \br^i,\bb\rangle^42\pi /n^2$. 
As $\mu_{\bU}$ is the conditional expectation of $M_{\bU}^{\sym}$, 
the Gaussian concentration inequality \citep{borell1975brunn,kwapien1994remark} provides 
\bes
\P_{\bU}\big\{M_{\bU}^{\sym} > \mu_{\bU} + \sigma_{\bU}t \big\} \le \P\big\{N(0,1)>t\big\},\ \forall\, t\ge 0. 
\ees
Thus, symmetrizing via the Jensen inequality as in \eqref{Dudley-0}, we have 
\bes
\E\big[e^{M_{\bU}/\eta}\big] 
\leq \E\big[\exp\big(M_{\bU}^{\sym}/\eta\big)\big] 
\leq \E\big[e^{\mu_{\bU}/\eta}\big(1/2+ \exp( \sigma_{\bU}^2/(2\eta^2))\big)\big]. 
\ees
As $\|\bb\|_1\le s_1^{1/2}$ for $\bb\in\scrC_1^*(s_1)$, 
we have $\sigma_{\bU}^2 \le M_{\bU}(2\pi s_1L_0/n)$. 
Thus, for $2\pi s_1L_0/n \le \eta/2$ and $\eta\le 1$, the above inequality and \eqref{new-pf-1} yield 
\bes
\E\big[e^{M_{\bU}/\eta}\big] 
&\leq& \E\big[e^{M_{\bU}/(4\eta)+(1+\eta)/4}\big(1/2+ e^{M_{\bU}/(4\eta)}\big)\big]
\cr &\leq& \big(\E\big[e^{M_{\bU}/\eta}\big]\big)^{1/4}(\sqrt{e}/2) 
+ \big(\E\big[e^{M_{\bU}/\eta}\big]\big)^{1/2}\sqrt{e}. 
\ees
For $x =\big(\E\big[e^{M_{\bU}/\eta}\big]\big)^{1/2}$, this gives $x \le e^{1/2}(1/(2\sqrt{x})+1)$, which implies $x\le \sqrt{5}$.  
Consequently, for $4\pi s_1L_0/n\le \eta\le 1$, $\P\big\{M_{\bU} > c_0 \big\} 
\le e^{ - c_0/\eta}\E e^{M_{\bU}/\eta}\le 5e^{-c_0/\eta}$. 
The conclusion follows as $c_0\in (0,1)$ and $5e^{-c_0/\eta} \le \eps_1\le 1$ imply $\eta\le 1$. 
 $\hfill\square$


\medskip
{\sc Proof of Lemma \ref{lemma-1}.} 
Let $\bM_i = \br_k^{i}\otimes \br_k^{i} - \E[\br_k^{i} \otimes \br_k^{i} ]$ 
where $\bu\otimes \bv = \bu\bv^\top$ for all vectors $\bu$ and $\bv$. 
We have $\|\bM_i\|_S \leq L_0^2d_{k}$ and 
$\sum_{i=1}^n \E \bM_i^2/n \leq L_0^2 d_{k} \bV_{k}$ 
with $\bV_{k} = \E\big( \bU_{k}^\top\bU_{k}/n \big)$. 
By the non-commutative Bernstein Inequality \citep{tropp2012user}, 
$$
   \P\big\{ \big\| \bU_{k}^\top\bU_{k}/n - \E \bV_{k}\big\|_S > x\big\}
  \leq 2d_{k}
  \exp\bigg[ \frac{ - n x^2/2}{L_0^2d_{k}\|\bV_{k}\|_{S} + L_0^2d_{k}x/3\}}\bigg] 
$$
for all $x>0$. This gives \eqref{lm-1-1} via the union bound. 
For \eqref{lm-1-2}, we have $\|\bV_k^{-1/2}\br^i_k\|_2\le \nu_-^{-1/2}L_0d_k^{1/2}$. 
Applying the above inequality to $(\bU_{k}\bV_{k}^{-1/2})$, we find that 
$$
   \P\bigg\{\max_{\E[\|\bU_k\bb_k\|_{2,n}^2]> 0} 
   \bigg| \frac{\|\bU_k\bb_k\|_{2,n}^2}{\E[\|\bU_k\bb_k\|_{2,n}^2]} - 1\bigg|  > x\bigg\}
  \leq 2d_{k}
  \exp\bigg[ \frac{ - n x^2/2}{(L_0^2/\nu_-)d_{k}(1+ x/3)}\bigg] 
$$
as the above maximum is $\big\|(\bU_{k}\bV_{k}^{-1/2})^\top(\bU_{k}\bV_{k}^{-1/2})/n - \bI_{d_k}\big\|_S$. 
$\hfill\square$

\medskip
{\sc Proof of Theorem \ref{th-5}.} 

(i) In $\Omega_0$, \eqref{th-1-1} holds. In $\Omega_1$, 
$$
C_{\rm pred}(\xi,S_0) \le 1/\kappa^2(\xi,S_0) \le (1+c_0)/\{(1-c_0)\kappabar^2(\xi_0,S_0)\}. 
$$ 
By \eqref{Omega_2}, 
$\Omega_2 = \big\{\Deltabar_{S_0} \le \sigma^2s_2/n \big\}$ 
for the $\Deltabar_{S_0}$ in Theorem \ref{th-1}. Thus, \eqref{th-1-1} implies \eqref{th-5-1}. 
Here $\xi = (A_0+1)/(A_0-1)$ is as in Theorem \ref{th-1}. 

(ii) Consider the case where $\Omega_0\cap\Omega_1\cap\Omega_2$ happens. 
Let $\bffbar = \sum_{(j,k)\in\scrK^*}\bU_{j,k}\bbeta^*_{j,k}$, 
$\bh_{j,k}=\hbbeta_{j,k}-\bbeta^*_{j,k}$ and $\bh = (\bh_{j,k}^\top, (j,k)\in\scrK^*)^\top$ 
with the $\scrK^*$ in \eqref{nu_pm}.    
As in the proof of Theorem~\ref{th-1} (i), 
in the event $\Omega_0\cap\Omega_2$. 
\bes
&& \big\|\hbf - \bffbar\big\|_{2,n}^2 + \big\|\hbf - \bff^*\big\|_{2,n}^2 
+ 2(A_0-1)\textsum_{(j,k)\in \scrK^*}\lam_{j,k}\|\bU_{j,k}\bh_{j,k}\|_{2,n}
\cr &\le& \sigma^2s_2/n + 4A_0\textsum_{(j,k)\in S_0}\lam_{j,k}\|\bU_{j,k}\bh_{j,k}\|_{2,n}. 
\ees

We shall consider two cases. In the first case where 
\bel{pf-th-5-1}
\sigma^2s_2/n \le 4A_0\textsum_{(j,k)\in S_0}\lam_{j,k}\|\bU_{j,k}\bh_{j,k}\|_{2,n},
\eel
we have $2(A_0-1)\textsum_{(j,k)\in S_0^c}\lam_{j,k}\|\bU_{j,k}\bh_{j,k}\|_{2,n}\le 
(6A_0+2)\textsum_{(j,k)\in S_0}\lam_{j,k}\|\bU_{j,k}\bh_{j,k}\|_{2,n}$, 
so that the proof of Theorem \ref{th-1} proceeds verbatim with 
$\xi = (3A_0+1)/(A_0-1)$. In particular, $\bh\in \scrC(\xi,S_0)$ 
and by the first inequality in $\Omega_1$ in \eqref{pf-th5-key}, 
\bes
\|\fhat - \fbar\|_{L_2}^2=\|\bU\bh\|_{L_2,n}^2 \le \|\bU\bh\|_{2,n}^2/(1-c_0) 
= \|\hbf-\bffbar\|_{2,n}^2/(1-c_0).  
\ees
In the second case where \eqref{pf-th-5-1} does not happen, we have 
\bes
\big\|\bh\big\|_1 
&\le& \sum_{(j,k)\in \scrK^*}\frac{n^{1/2}\lam_{j,k}}{\sigma}\|\bh_{j,k}\|_{2}
\cr &\le&
 \sum_{(j,k)\in \scrK^*}\frac{n^{1/2}\lam_{j,k}}{\sigma \nu_-^{1/2}}\|\bU_{j,k}\bh_{j,k}\|_{L_2,n}
\cr &\le& \sum_{(j,k)\in \scrK^*}\frac{n^{1/2}\lam_{j,k}}{\sigma \nu_-^{1/2}}
\frac{\|\bU_{j,k}\bh_{j,k}\|_{2,n}}{1-c_0}
\cr &\le& \frac{\sigma n^{-1/2}s_2}{(A_0-1)\nu_-^{1/2}(1-c_0)}
\ees
by \eqref{lam}, \eqref{nu_pm} and the third inequality in $\Omega_1$ in \eqref{pf-th5-key}. 
Thus, by the second inequality in $\Omega_1$ 
with $\bb = (A_0-1)\nu_-^{1/2}(1-c_0)(n/s_2)^{1/2}\bh$ and under the condition $s_2\le s_1$ 
\bes
\|\fhat - \fbar\|_{L^2} 
= \begin{cases} \|\bU \bh\|_{L_2,n}^2 \le \|\bU \bh \|_{2,n}^2/(1-c_0), & \|\bU\bb\|_{L_2,n}^2 \ge 1,
\cr \|\bU \bh\|_{L_2,n}^2 < \sigma^2 s_2/\{(1-c_0)^2(A_0-1)^2\nu_-n\}, & \|\bU\bb\|_{L_2,n}^2 < 1.
\end{cases}
\ees
Thus, in either cases, \eqref{th-5-1} and the above bounds for $\|\fhat - \fbar\|_{L^2}$ yield
\bes
&& \|\fhat - \fbar\|_{L_2}^2 + \|\fhat - f^*\|_{L_2}^2 
\cr &\le& 3\|\fhat - \fbar\|_{L_2}^2 + 2\|\fbar-f^*\|_{L_2}
\cr &\le& \max\bigg(\frac{3\sigma^2 s_2}{(1-c_0)^2(A_0-1)^2\nu_-n}, 
\frac{3\|\hbf-\bffbar\|_{2,n}^2}{1-c_0}\Big) + \frac{\sigma^2 s_2}{n}, 
\cr &\le& C_{A,\nu_-,c_0}\big(\|\lam_S\|_2^2/\kappabar^2(\xi_0,S_0)+\sigma^2s_2/n\big), 
\ees
with a constant $C_{A_0,\nu_-,c_0}$ depending on $(A_0,\nu_-,c_0)$ only. 
This gives \eqref{th-5-2}. 

(iii) By \eqref{lam-gen} and Theorem~\ref{th-1}~(ii), we have $\P\{\Omega_0^c\} \le \eps/\sqrt{2\log(p/\eps)}$. 
By \eqref{th-4-2} and \eqref{lm-1-2} we have 
$\P\{\Omega_1^c\} \le 2\eps_1$ for the $\Omega_1$ 
in \eqref{pf-th5-key} when \eqref{cond-CC} hold with $\{c_0,\eps_1\}\subset(0,1)$ and $s_2\le s_1$. 
It remains to prove $\P\{\Omega_1\setminus \Omega_2\} \le \eps_2$
for the $s_2$ in \eqref{s_2}. 
Because 
$$
\sum_{(j,k)\in\scrK^*\setminus S_0} \lam_{j,k}\|\bU_{j,k}\bbeta^*_{j,k}\big\|_{2,n} \le 
(1+c_0)\nu_+^{1/2} \sum_{(j,k)\in \scrK^*\setminus S_0}\lam_{j,k}\|\bbeta^*_{j,k}\|_2
$$ 
by \eqref{nu_pm} and the third inequality in $\Omega_1$ in \eqref{pf-th5-key} 
and $\P\big\{ \big\|\bff^* - \bffbar\big\|_{2,n}^2 \ge \| f^*-\fbar\|_{L_2}/\eps_2\big\}\le\eps_2$ 
by the Markov inequality, it suffices to prove 
\bel{pf-th-5-2}
\P\bigg\{ \big\|\bff^* - \bffbar\big\|_{2,n}^2 \ge  
2\| f^*-\fbar\|_{L_2} 
+ \frac{(C_{\alpha-1/2}L_0M_{\alpha,1})^2|\log \eps_2|}{n^{2(\alpha+\alpha^*)/(2\alpha_*+1)}}
\bigg\} \le \eps_2. 
\eel
in view of \eqref{Omega_2} and \eqref{s_2}. 
As $\big\|\bff^* - \bffbar\big\|_{2,n}^2$ is an average of independent variables 
$(f^*(\br^i)-\fbar(\br^i))^2$ with $\br^i$ being the rows of $\bX$, the Bernstein inequality yields 
\bel{pf-th-5-3}
\P\{\big\|\bff^* - \bffbar\big\|_{2,n}^2 
\ge  \mu_* + \sqrt{\mu_*C_*2|\log \eps_2|/n} + 2C_*|\log \eps_2|/(3n) \Big\} 
\le \eps_2
\eel
with $\mu_* = \E\big[\|\bff^* - \bffbar\|_{2,n}^2\big] = \|f^*- \fbar\|_{L_2}^2$ and 
$C_* = \|f^* - \fbar\|_\infty^2$. Here the variance bound 
$n\Var\big(\|\bff^* - \bffbar\|_{2,n}^2\big)
\le \mu_*C_*$ is used. 
 
As $\|u_{j,k,\ell}(\cdot)\|_\infty\le L_0$ and $f^*-\fbar=\sum_{j=1}^p\sum_{k>k^*}\sum_{\ell=1}^{2^{k-1}}u_{j,k,\ell}(x_j)\beta^*_{j,k,\ell}$, 
\eqref{M_alpha} gives 
$$
C_*
\le \bigg(\sum_{j=1}^p \sum_{k > k^*} 2^{(k-1)/2}L_0\|\bbeta^*_{j,k}\|_2\bigg)^2
\le \frac{1}{2}\bigg(\sum_{j=1}^p \frac{C_{\alpha-1/2}L_0\|f^*\|_{\alpha,2}}{2^{(\alpha-1/2)k^*/2}}\bigg)^2
\le \frac{(C_{\alpha-1/2}L_0M_{\alpha,1})^2}{2n^{2(\alpha-1/2)/(2\alpha_*+1)}}
$$ 
with $C_\alpha = 1/(4^\alpha-1)^{1/2}$, in view of the condition $2^{k^*}\ge n^{1/(2\alpha^*+1)}$.  
Thus, 
\bes
\mu_* + \sigma_*\sqrt{2|\log \eps_2|/n} + 2C_*|\log \eps_2|/(3n)
\le 2\mu_*+ \bigg(\frac{1}{2}+\frac{1}{3}\bigg)\frac{C_*|2\log \eps_2|}{n}
\ees
by the bounds for $\mu_*$ and $C_*$. This and \eqref{pf-th-5-3} yield \eqref{pf-th-5-2}. 

(iv) We notice that when $L_0$ and $c_0$ are fixed positive constants 
and both $1/\nu_-$ and $1/\kappabar (\xi_0,S_0)$ are bounded, 
we are allowed to take $s_1\le \max\{O(s),s_2\}$ and for this $s_1$ 
\eqref{cond-CC} holds with $\eps_1=o(1)$ when $n \gg 2^{k^*}\log(np)$ 
and $n \gg s_1(\log s_1)^2(\log n)\log(np)$. Thus, part (iv) follows from parts (ii) and (iii).  

(v) We omit the proof of this part as it is identical to the proof of Theorem 2. $\hfill\square$

\bibliographystyle{apalike}
\bibliography{MR-Lasso}

\end{document}

%% file: 04-pre-am-short.tex
\def\argmin{\mathop{\rm arg\, min}}

\newcommand{\bel}{\begin{eqnarray}\label}
\newcommand{\eel}{\end{eqnarray}}
\newcommand{\bes}{\begin{eqnarray*}}
\newcommand{\ees}{\end{eqnarray*}}
\newcommand{\bei}{\begin{itemize}}
\newcommand{\eei}{\end{itemize}}
\newcommand{\beiftnt}{\begin{itemize}\footnotesize}

\def\benu{\begin{enumerate}}
\def\eenu{\end{enumerate}}

\def\argmin{\mathop{\rm arg\, min}}
\def\textsum{\hbox{$\sum$}}

\def\real{{\mathbb{R}}}
\def\R{{\real}}

\def\E{{\mathbb{E}}}

\def\P{{\mathbb{P}}}

\def\complex{\mathop{{\rm I}\kern-.58em\hbox{\rm C}}\nolimits}

\def\rank{\hbox{\rm rank}}

\def\Var{\hbox{\rm Var}}

\def\pen{\hbox{\rm pen}}

\def\mathbold{\boldsymbol} 

\def\ba{\mathbold{a}}

\def\bb{\mathbold{b}}

\def\calB{{\cal B}}

\def\scrC{{\mathscr C}}

\def\scrC{{\mathscr C}}

\def\bff{\mathbold{f}}\def\fhat{\widehat{f}}
\def\hbf{{\widehat{\bff}}}
\def\scrF{{\mathscr F}}

\def\bg{\mathbold{g}}

\def\bh{\mathbold{h}}

\def\bI{\mathbold{I}}

\def\kbar{{\overline k}}

\def\scrK{{\mathscr K}}

\def\scrL{{\mathscr L}}

\def\bM{\mathbold{M}}

\def\bP{\mathbold{P}}

\def\br{\mathbold{r}}

\def\tbar{{\overline t}}

\def\bu{\mathbold{u}}

\def\bU{\mathbold{U}}

\def\bv{\mathbold{v}}

\def\bV{\mathbold{V}}

\def\bx{\mathbold{x}}

\def\bX{\mathbold{X}}

\def\by{\mathbold{y}}

\def\bbeta{\mathbold{\beta}}\def\hbeta{\widehat{\beta}}
\def\betabar{{\overline{\beta}}}
\def\hbbeta{{\widehat{\bbeta}}}
\def\bbetabar{{\overline\bbeta}}

\def\Deltabar{{\overline \Delta}}

\def\ep{\varepsilon}\def\eps{\epsilon}\def\veps{\varepsilon}
\def\bep{\mathbold{\ep}}

\def\kappabar{{\overline{\kappa}}}

\def\lam{\lambda}